\documentclass[11pt, reqno]{amsart}
\usepackage{pdfpages}
\usepackage{amsmath}
\usepackage{fullpage}
\usepackage{amssymb}
\usepackage{amsthm}
\usepackage{bbm}
\usepackage{subfigure}
\usepackage{float}
\usepackage{graphicx}
\usepackage{color}

\usepackage{hyperref}
\usepackage{enumerate}
\graphicspath{{../figs/}}
\newtheorem{lemma}{Lemma}[section]
\newtheorem{theorem}[lemma]{Theorem}
\newtheorem{proposition}[lemma]{Proposition}
\newtheorem{remark}[lemma]{Remark}

\DeclareMathOperator{\re}{Re}
\DeclareMathOperator{\im}{Im}
\numberwithin{equation}{section}
\newcommand{\fe}{\mathrm{e}}

\title{Error estimates of time-splitting schemes for nonlinear Klein--Gordon equation with rough data}

\author{Lun Ji}
\address{LSEC, ICMSEC, Academy of Mathematics and Systems Science,
Chinese Academy of Sciences, Beijing 100049, China (L.~Ji)}
\address{School of Mathematical Sciences, University of Chinese Academy of Sciences, Beijing 100049, China (L.~Ji)}
\email{ujeybn@lsec.cc.ac.cn}

\author{Xiaofei Zhao}
\address{School of Mathematics and Statistics \& Computational Sciences Hubei Key Laboratory, Wuhan University, 430072 Wuhan, China (X.~Zhao)}
\email{matzhxf@whu.edu.cn}

\begin{document}

\subjclass[2010]{65M12, 65M15, 65M70, 81Q05, 35Q40.}

\keywords{Nonlinear Klein--Gordon equation, rough solution, discrete Bourgain space, time-splitting scheme, optimal error estimate}

\date{}
\dedicatory{}

\begin{abstract}
In this work, we consider the convergence analysis of time-splitting schemes for the nonlinear Klein--Gordon/wave equation under rough initial data. The optimal error bounds of the Lie splitting and the Strang splitting are established with sharp dependence on the regularity index of the solution from a wide range that is approaching the lower bound for well-posedness. Particularly for very rough data, the technique of discrete Bourgain space is utilized and developed, which can apply for general second-order wave models. Numerical verifications are provided. 
\end{abstract}

\maketitle

%%%%%%%%%%%%%%%%%%%%%%%%%%%
\section{Introduction}
%%%%%%%%%%%%%%%%%%%%%%%%%%%

The nonlinear Klein--Gordon/wave equation is mathematically one of the most important dispersive models \cite{Coweb,Taoams}
\begin{equation}\label{kg}
\begin{aligned}
&z_{tt}=\Delta z-m\,z+\lambda z^3, \quad z=z(t,x),\quad t>0,\ x\in\mathbb{T}^d,\\
&z(0,x)=z_0(x),\quad z_t(0,x)=z_1(x),\quad x\in\mathbb{T}^d,
\end{aligned}
\end{equation}
where $z(t,x)$ is the real-valued unknown scalar field with $z_0,z_1$ two given initial data, $()_t$ denotes $\partial_t()$, $\mathbb{T}^d$ denotes a $d$-dimensional torus with $d=1,2,3$, and $ m\geq0,~\lambda\in\mathbb{R}$ are given parameters. It is widely used in quantum physics, plasma physics, cosmology, and general wave motions \cite{BergeColin,Davydov,Huang,Mauser}. The torus domain indicates periodic boundary conditions for $z$, which also serves as a valid truncation of the whole space problem when the initial localized waves are yet to reach the boundary of a computational domain. When $m>0$ and $\lambda\neq0$, (\ref{kg}) is then the nonlinear Klein--Gordon equation (NKG) with $\lambda>0$ and $\lambda<0$ representing the focusing and defocusing nonlinear interaction, respectively, which can be interpreted as the nonlinear Schr\"odinger equation (NLS) integrated with relativistic effect \cite{BaomZhao}. The energy is conserved by (\ref{kg}) in time as
$$
E(t):=\int_{\mathbb{T}^d}\left[|z_t|^2+|\nabla z|^2+m|z|^2-\frac{\lambda}{2}|z|^4\right]dx\equiv E(0),\ t\geq0.
$$

Theoretically, \eqref{kg} has been intensively studied. We refer to \cite{Coweb, Charxiv,Lindjfa,Taoams} for its local well-posedness theorems and \cite{Bourgain,Coweb,Ginibre1,Ginibre2,Taoams} for global ones, where the results strongly depend on the regularity $s>0$ of the initial data $(z_0,z_1)\in H^s\times H^{s-1}$ and the dimension $d$. Seeking the minimum requirement on $s$ for well-posedness has always been a challenging but meaningful goal, which enables the model to allow more realistic data with noise/roughness. For example, in two dimensions, i.e., $d=2$, $s>\frac14$ is proved as a sufficient condition \cite{Lindjfa} for local well-posedness and $s\geq\frac16$ is necessary, while the situation for $s\in(\frac16,\frac14]$ is still an open question. Note that since the speed of the propagation of $z$ is limited to the speed of light, the wave-type equations on the torus share the same well-posedness results with the whole space case \cite{Coweb}.

Numerically, \eqref{kg} has also been intensively solved using various types of numerical techniques. Early research interests focus on structure-preserving and high-order accuracy \cite{DG,Duncan,LiGuo,LiuWu}. However, when the solution is not smooth enough in space, popular traditional temporal discretizations, including the finite difference, exponential integrators, and time-splitting methods, will all suffer from convergence order reduction in time \cite{low} to some extent. Therefore, on the one hand, some low-regularity integrators have been proposed in recent years \cite{BDHoch,CaoLLY,Li,WangZhao}, aiming to accurately integrate \eqref{kg} under rough data. 
See such a kind of works on other dispersive equations \cite{FMW,Jimcom,LiWunls,LiWukdv,LRIDirac,WuYao}. On the other hand, it is also important to clearly understand the precise convergence result and the sharp accuracy order of the traditional discretizations for \eqref{kg} under rough data, which have not been done so far. Their rigorous numerical analysis, particularly when the regularity of the solution is very low, so that standard Sobolev embedding fails, is quite challenging.

Recently, discrete Strichartz estimates \cite{LiWukdv,Ostfocm} and Bourgain spaces \cite{Ostjems,Roufocm} have been established to numerically solve dispersive PDEs under extremely low-regularity constraints. Those tools are first established for one-dimensional NLS, then developed to higher dimensional cases \cite{Jisiam,Jiima} and other dispersive equations \cite{Jimcom,Roupaa}. However, unlike the Sobolev setting, the structure of these tools depends strongly on the linear operator of the equation, and for wave-type equations, such low-regularity estimates are still open.

This work is devoted to rigorously analyzing the convergence of time-splitting methods, which are one class of the most popular time integration techniques, to solve \eqref{kg} with rough initial data approaching the critical (minimal requirement for well-posedness) regularity index $s$. The time-splitting methods concerned here include the Lie-Trotter splitting scheme and the Strang splitting scheme, where the latter, in fact, also equivalently reads as a time-symmetric exponential integrator \cite{Dongcicp,Lubich}. To precisely describe the convergence result, we will mainly work on the two-dimensional case of \eqref{kg} for presentations, while the generalization to the other cases can be established in the same manner, which we will briefly discuss. Moreover, we will use discrete Bourgain spaces, which, as an alternative to Strichartz estimates in the wave setting, provide a convenient framework for the analysis of approaching critical regularity under periodic boundary conditions. Indeed, for NKG \eqref{kg}, in order to obtain the complete and optimal convergence order of the splitting schemes, we will work separately in the Sobolev space and in the Bourgain space, based on the range of $s$. In two dimensions, that is $s>\frac12$ and $s\in(\frac{11}{40},\frac{51}{40}]$. Particularly, some discrete wave-type Bourgain estimates are established, which would be useful also for more general wave models. Numerical experiments will be presented in the end to verify the sharpness of the theoretical estimates.

The rest of the paper is organized as follows. Section \ref{sec:scheme} reviews the splitting schemes for the NKG and Section \ref{sec:sobolev} presents the convergence analysis in Sobolev space. In Section \ref{sectbourg}, we present the analysis in the discrete Bourgain space. Some technical estimates and generalizations are given in Section \ref{sectproof}, and numerical illustrations are given in Section \ref{sectnumexp}. 

\subsection*{Notations}
We denote the time step size by $\tau$, and we use the Japanese bracket $\langle\cdot\rangle=(m+|\cdot|^2)^{\frac12}$ notation. The estimate $\lesssim$ means $A\leq CB$, where $C$ is a generic constant independent of $\tau$. The notation $\lesssim_\gamma$ emphasizes that the constant $C$ is particularly dependent on $\gamma$. Moreover, $A\sim B$ means that $A\lesssim B\lesssim A$. A basic embedding inequality reads 
\begin{equation}\label{cbe}
\|uv\|_{H^r}\lesssim\|u\|_{H^r}\|v\|_{H^s},\quad s>\tfrac d2,~0\leq r\leq s.
\end{equation}

For a function $f(x)=\sum\limits_{k\in\mathbb{Z}^d}\widehat{f_k}{\rm e}^{i\langle k, x\rangle}$, where $\langle\cdot,\cdot\rangle$ denotes the inner product in $\mathbb R^d$, the pseudo-differentiation operator is defined through Fourier space:
$
\langle\nabla\rangle^\alpha f=\sum\limits_{k\neq 0}\langle k\rangle^\alpha\widehat{f_k}{\rm e}^{i\langle k,x\rangle},\ \alpha\in\mathbb{R}.$
For a space-time function, we may omit the spatial variable for simplicity, e.g., $z(t)=z(t,x)$.

Finally, for sequences $(u_n)_{n\in\mathbb{Z}} \in X^\mathbb{Z}$, where $X$ is a Banach space equipped with norm $\|\cdot\|_{X}$, we employ the usual norms
$$
\|u_n\|_{l_\tau^pX}=\Big(\tau\sum_n\|u_n\|_X^p \Big)^\frac1p, \quad
\|u_n\|_{l_\tau^\infty X} = \sup_{n\in\mathbb{Z}} \|u_n\|_X.
$$

%%%%%%%%%%%%%%%%%%%%%%
\section{Time-splitting methods}\label{sec:scheme}
%%%%%%%%%%%%%%%%%%%%%%

Let us briefly review the time-splitting schemes for \eqref{kg} in this section. 
First of all, we can rewrite \eqref{kg} as
\begin{equation}\label{kga}
\left(\begin{array}{c}
z \\ z_t
\end{array}\right)_t=
\left(\begin{array}{cc}
0 & 1 \\ \Delta-m & 0
\end{array}\right)
\left(\begin{array}{c}
z \\ z_t
\end{array}\right)+\lambda
\left(\begin{array}{c}
0 \\ z^3
\end{array}\right),\quad
\left(\begin{array}{c}
z \\ z_t
\end{array}\right)(0,x)=
\left(\begin{array}{c}
z_0(x) \\ z_1(x)
\end{array}\right),
\end{equation}
and by left multiplying $(1,i\langle\nabla\rangle^{-1})$, we obtain \begin{equation}\label{kgb}
u_t=-i\langle\nabla\rangle u+i\lambda\langle\nabla\rangle^{-1}(\re u)^3,\quad u(0,x)=u_0(x)=z_0(x)+i\langle\nabla\rangle^{-1}z_1(x).
\end{equation}
where $u=z+i\langle\nabla\rangle^{-1}z_t$. Therefore, we can split \eqref{kgb} into two subsystems \cite{Baomcom,Dongcicp}:
\begin{align}
\label{kgsub1}v_t&=-i\langle\nabla\rangle v,\\
\label{kgsub2}w_t&=i\lambda\langle\nabla\rangle^{-1}(\re w)^3.
\end{align}
Here, we can see easily that the first subsystem \eqref{kgsub1} can be solved exactly as $v(t)=\fe^{-it\langle\nabla\rangle}v_0$. For the second subsystem \eqref{kgsub2}, note that $\langle\nabla\rangle^{-1}(\re w)^3\in\mathbb{R}$, we actually have $\re w_t=0$, i.e., $\re w=\re w_0$. Thus, \eqref{kgsub2} can also be solved exactly as $w(t)=w_0+i\lambda t\langle\nabla\rangle^{-1}(\re w_0)^3$. 
By composing the two exact solutions, we can get for \eqref{kgb} the Lie splitting scheme
\begin{equation}\label{kglie}
u_{n+1}=\Phi(u_n)=\fe^{-i\tau\langle\nabla\rangle}u_n+i\lambda\tau\langle\nabla\rangle^{-1}\fe^{-i\tau\langle\nabla\rangle}(\re u_n)^3,
\end{equation}
and the Strang splitting scheme
\begin{equation}\label{kgstrang}
u_{n+1}=\Psi(u_n)=\fe^{-i\tau\langle\nabla\rangle}u_n+i\lambda\tau\langle\nabla\rangle^{-1}\fe^{-i\tau\langle\nabla\rangle/2}(\re \fe^{-i\tau\langle\nabla\rangle/2}u_n)^3,
\end{equation}
where $\tau>0$ denotes the time step size. Equivalently, the Strang splitting scheme (\ref{kgstrang}) can be derived based on the Deuflhard-type exponential integration \cite{Dongcicp}.

For the error estimate work in the rest of the paper, without loss of generality, we will simply fix the parameters in (\ref{kg}) as $m=1,\ \lambda=-1$, i.e., the defocusing NKG. The case where $m=0$ is briefly discussed in Section~\ref{remhighd}. We will work on the two-dimensional case, i.e., $d=2$ in (\ref{kg}), and provide remarks about the one-dimensional and three-dimensional cases. The two-dimensional defocusing NKG, as known from the literature \cite{Lindjfa}, is locally well-posed for $(z_0,z_1)\in H^s\times H^{s-1}$ with $s>\frac14$.
Our study will therefore be up to a fixed finite time $T>0$ within the maximum lifespan of the solution. 
We shall carry out the analysis in the next two sections for $s>\frac12$ in the Sobolev space and for $\frac{11}{40}<s\leq\frac{51}{40}$ in Bourgain space. 

%%%%%%%%%%%%%%%%%%%%%%
\section{Error estimates in Sobolev space}\label{sec:sobolev}
%%%%%%%%%%%%%%%%%%%%%%

In this section, we will study the convergence behavior of the splitting methods \eqref{kglie} and \eqref{kgstrang} when the solution of the NKG \eqref{kg} is not too rough. The main efforts are devoted to establishing the following theorem, which will be done by analysis in the subsections in a sequel.
\begin{theorem}[Convergence in Sobolev]\label{theosob}
For $s>\frac12$ and the initial data $z_0(x)\in H^s(\mathbb{T}^2),~z_1(x)\in H^{s-1}(\mathbb{T}^2)$, let $z$ be the exact solution to \eqref{kg} on $[0,T]$. Moreover, we denote $u=z+i\langle\nabla\rangle^{-1}z_t$ (i.e., $u$ is the solution to the equivalent equation \eqref{kgb}). Furthermore, let $u_n,~v_n$ be the sequences defined by the Lie splitting \eqref{kglie} and the Strang splitting \eqref{kgstrang}, respectively. Then, there exist some constants $\tau_0,~C_T>0$ such that for any time step size $\tau\in(0,\tau_0]$ and $0\leq n\tau\leq T$,
\begin{subequations}
\begin{align}
\label{liesob}\mbox{Lie:  }\quad\|u(t_n)-u_n\|_{H^\frac12}&\leq C_T\tau^{3s-\frac32},\quad s\in(\tfrac12,\tfrac56],\\
\label{strangsob}\mbox{Strang:  }\quad\|u(t_n)-v_n\|_{H^\frac12}&\leq C_T\tau^{s+\frac12},\quad s\in(1,\tfrac32],
\end{align}
\end{subequations}
where $C_T$ depends on $T$ and $\tau_0$ but is independent of $n$ and $\tau$.
\end{theorem}

\subsection{Cauchy problem for (\ref{kgb})}

We begin with the discussion and analysis of the Cauchy problem for \eqref{kgb}, which will serve as needed prior estimates for error estimates later. Technically, the classical bilinear estimate \eqref{cbe} is not enough for rigorous analysis here. Therefore, we give the following tool lemma, which is an extension of \eqref{cbe}.
\begin{lemma}\label{theomultsob}
For any functions $u,v,w$ defined in $H^{s_i}(\mathbb{T}^2),~i=1,2,3$, respectively, we have the following estimate:
\begin{equation}\label{multsob}
\|uvw\|_{H^r}\lesssim\|u\|_{H^{s_1}}\|v\|_{H^{s_2}}\|w\|_{H^{s_3}},\quad r\in(-1,1),~s_i\in(r,1),~s_1+s_2+s_3\geq r+2.
\end{equation}
\end{lemma}
We postpone the proof of this lemma to Section~\ref{sectprosob}. By using it, 
we can now establish the existence and uniqueness of the exact solution to \eqref{kgb}.

\begin{proposition}\label{theoexist}
Let $u_0\in H^s(\mathbb{T}^2)$ with $s\geq\frac12$. Then, there exist some $T^*\in(0,+\infty]$ and a unique maximal solution $u\in\mathcal{C}([0,T],H^s(\mathbb{T}^2))$ satisfying \eqref{kgb}, for every $T\in(0,T^*)$.
\end{proposition}
\begin{proof}
Consider the Duhamel form of \eqref{kgb}: 
$
v= F(v),
$
where
\begin{equation}\label{duh}
F(v)(t)=\fe^{-it\langle\nabla\rangle}u_0-i\langle\nabla\rangle^{-1}\int_0^t\fe^{-i(t-\vartheta)\langle\nabla\rangle}(\re v(\vartheta))^3d\vartheta.
\end{equation}
Thus by using \eqref{cbe} and \eqref{multsob}, we deduce that
$$
\|F(v)\|_{Y^s(t)}\leq \|u_0\|_{H^s(\mathbb{T}^2)}+t\|v^3\|_{Y^{s-1}(t)}\leq\|u_0\|_{H^s(\mathbb{T}^2)}+Ct\|v\|^3_{Y^s(t)},
$$
where we denote $Y^s(t)=L^\infty H^s([0,t]\times\mathbb{T}^2)$ for short, and $C>0$ is a constant that only depends on $s$. Moreover, from the same argument, we can deduce
\begin{equation}\label{compress}
\|F(v_1)-F(v_2)\|_{Y^s(t)}\leq 3CtR^2\|v_1-v_2\|_{Y^s(t)},
\end{equation}
for all $v_1,v_2$ satisfying $\|v_i\|_{Y^s(t)}\leq R$, where $R>\|u_0\|_{H^s}$. Thus by the Banach fixed-point theorem, we can get the existence of a fixed point of $F$ (i.e., a solution of \eqref{kgb}) in the ball of radius $R$ in $Y^s(t)$ for $t\in(0,\frac1{3CR^2})$.

By the same argument, we can also get the uniqueness of solution in $Y^s(t)$ for $s\geq\frac12$. Moreover, by a standard iteration, we can get a unique solution to \eqref{kgb} in $Y^s(T)$, then the maximal solution follows in a standard way.
\end{proof}

We then move on to the error estimates.

\subsection{Local error estimate}\label{sectsob}

In this part, we estimate the local error of the Lie splitting \eqref{kglie} and the Strang splitting \eqref{kgstrang}. First of all, by \eqref{kglie}, \eqref{kgstrang} and \eqref{duh}, the local error of Lie and Strang can be written respectively as
\begin{equation}\label{clieloc}
\begin{aligned}
\mathcal{E}_t^1(u(t_n))&:=\Phi(u(t_n))-u(t_{n+1})\\
&=-i\tau\langle\nabla\rangle^{-1}\fe^{-i\tau\langle\nabla\rangle}(\re u(t_n))^3+i\langle\nabla\rangle^{-1}\int_0^\tau\fe^{-i(\tau-\xi)\langle\nabla\rangle}(\re u(t_n+\xi))^3d\xi\\
&=i\langle\nabla\rangle^{-1}\fe^{-i\tau\langle\nabla\rangle}\int_0^\tau(\fe^{i\xi\langle\nabla\rangle}-1)(\re u(t_n+\xi))^3d\xi\\
&\quad+i\langle\nabla\rangle^{-1}\fe^{-i\tau\langle\nabla\rangle}\int_0^\tau\left[(\re u(t_n+\xi))^3-(\re u(t_n))^3\right]d\xi\\
&=:I_1(t_n)+I_2(t_n),
\end{aligned}
\end{equation}
\begin{align}
\mathcal{E}_t^2(u(t_n))&:=\Psi(u(t_n))-u(t_{n+1})\nonumber\\
&=-i\tau\langle\nabla\rangle^{-1}\fe^{-i\tau\langle\nabla\rangle/2}(\re \fe^{-i\tau\langle\nabla\rangle/2}u(t_n))^3\nonumber\\
&\quad+i\langle\nabla\rangle^{-1}\int_0^\tau\fe^{-i(\tau-\xi)\langle\nabla\rangle}(\re u(t_n+\xi))^3d\xi\nonumber\\
&=i\tau\langle\nabla\rangle^{-1}\fe^{-i\tau\langle\nabla\rangle/2}\left[(\re u(t_n+\tfrac\tau2))^3-(\re \fe^{-i\tau\langle\nabla\rangle/2}u(t_n))^3\right]\nonumber\\
&\quad+i\langle\nabla\rangle^{-1}\fe^{-i\tau\langle\nabla\rangle/2}\int_0^\tau \left[\fe^{-i(\tau/2-\xi)\langle\nabla\rangle}(\re u(t_n+\xi))^3-(\re u(t_n+\tfrac\tau2))^3\right]d\xi\nonumber\\
&=:I_3(t_n)+I_4(t_n).\label{cstrangloc}
\end{align}

We then have the following proposition for temporal local errors.
\begin{proposition}
For $s>\frac12$, with $u$ from Proposition~\ref{theoexist} and $\Phi,~\Psi$ defined in \eqref{kglie} and \eqref{kgstrang}, we have the following local error estimates:
\begin{subequations}
\begin{align}
\label{lieloc1}\|\mathcal{E}_t^1(u(t_n))\|_{H^\frac12}&\leq C_T\tau^{3s-\frac12},\quad s\in(\tfrac12,\tfrac56],\\
\label{strangloc1}\|\mathcal{E}_t^2(u(t_n))\|_{H^\frac12}&\leq C_T\tau^{s+\frac32},\quad s\in(1,\tfrac32].
\end{align}
\end{subequations}
\end{proposition}

\begin{proof}
First of all, by Proposition~\ref{theoexist}, we have
\begin{equation}\label{boundu}
\sup\limits_{\xi\in[0,\tau]}\|u(t_n+\xi)\|_{H^s}\leq \|u\|_{L^\infty H^s([0,T]\times\mathbb{T}^2)}\leq C_T,
\end{equation}
where $C_T$ is a constant that depends on $s,~\|u_0\|_{H^s}$ and $T$. For $s\in(\frac12,\frac56]$, take $r=3s-2$, then we have $r\in(-\frac12,\frac12]$. Thus by \eqref{multsob} and \eqref{boundu}, we have
\begin{equation}\label{e11}
\begin{aligned}
\|I_1(t_n)\|_{H^\frac12}&\leq\tau\sup\limits_{\xi\in[0,\tau]}\|(\fe^{i\xi\langle\nabla\rangle}-1)(\re u(t_n+\xi))^3\|_{H^{-\frac12}}\\
&\leq\tau\sup\limits_{\xi\in[0,\tau]}\|(\xi\langle\nabla\rangle)^{r+\frac12}(u(t_n+\xi))^3\|_{H^{-\frac12}}\\
&\lesssim\tau^{r+\frac32}\sup\limits_{\xi\in[0,\tau]}\|u(t_n+\xi)^3\|_{H^r}\lesssim\tau^{3s-\frac12}\sup\limits_{\xi\in[0,\tau]}\|u(t_n+\xi)\|^3_{H^s}\lesssim\tau^{3s-\frac12}.
\end{aligned}
\end{equation}

Moreover, by \eqref{multsob}, \eqref{duh} and \eqref{boundu}, we deduce that
\begin{equation}\label{e12}
\begin{aligned}
\|I_2(t_n)\|_{H^\frac12}&\leq\tau\sup\limits_{\xi\in[0,\tau]}\|(u(t_n+\xi))^3-(u(t_n))^3\|_{H^{-\frac12}}\\
&\lesssim\tau\sup\limits_{\xi\in[0,\tau]}\|u(t_n+\xi)-u(t_n)\|_{H^{\frac32-2s}}\sup\limits_{\xi\in[0,\tau]}\|u(t_n+\xi)\|^2_{H^s}\\
&\lesssim\tau\sup\limits_{\xi\in[0,\tau]}\|(\fe^{-i\xi\langle\nabla\rangle}-1)u(t_n)\|_{H^{\frac32-2s}}+\tau^2\sup\limits_{\xi\in[0,\tau]}\|u(t_n+\xi)^3\|_{H^{-\frac12}}\\
&\lesssim\tau^{1+3s-\frac32}\|u(t_n)\|_{H^s}+\tau^2\sup\limits_{\xi\in[0,\tau]}\|u(t_n+\xi)\|^3_{H^\frac12}\lesssim\tau^{3s-\frac12}.
\end{aligned}
\end{equation}
Combining \eqref{e11} and \eqref{e12}, we obtain \eqref{lieloc1}.

For the Strang splitting, when $s\in(1,\frac32]$, by \eqref{cbe}, \eqref{multsob}, \eqref{duh}, we have
\begin{equation}\label{e13}
\begin{aligned}
\|I_3(t_n)\|_{H^\frac12}&\leq\tau\|(\re u(t_n+\tfrac\tau2))^3-(\re \fe^{-i\tau\langle\nabla\rangle/2}u(t_n))^3\|_{H^{-\frac12}}\\
&\lesssim\tau\|\!\re\big(u(t_n+\tfrac\tau2)-\fe^{-i\tau\langle\nabla\rangle/2}u(t_n)\big)\|_{H^{-\frac12+2\varepsilon}}\sup\limits_{\xi\in[0,\tau]}\|u(t_n+\xi)\|^2_{H^{1-\varepsilon}}\\
&\lesssim\tau\big\|\!\im\big(\textstyle\int_0^\frac\tau2 \fe^{-i(\tau/2-\xi)\langle\nabla\rangle}\re (u(t_n+\xi))^3d\xi-\int_0^\frac\tau2\re(u(t_n+\xi))^3d\xi\big)\big\|_{H^{-1+\varepsilon}}\\
&\lesssim\tau^2\sup\limits_{\xi\in[0,\tau]}\|(\fe^{-i(\tau/2-\xi)\langle\nabla\rangle}-1)\re (u(t_n+\xi))^3\|_{H^{-1+\varepsilon}}\\
&\lesssim\tau^3\sup\limits_{\xi\in[0,\tau]}\|u(t_n+\xi)^3\|_{H^\varepsilon}\lesssim\tau^3\sup\limits_{\xi\in[0,\tau]}\|u(t_n+\xi)\|^3_{H^s}\lesssim\tau^3,
\end{aligned}
\end{equation}
where $\varepsilon\in(0,\frac12)$ can be taken arbitrarily small. Note that in the third line, we have freely added the integral $\int_0^\frac\tau2(\re u(t_n+\xi)^3)d\xi$, since $\im(\re u(t_n+\xi)^3)=0$.

If we denote $f(\xi)=\langle\nabla\rangle^{-1}(\fe^{-i(\tau/2-\xi)\langle\nabla\rangle}(\re u(t_n+\xi))^3)$, then $f$ is continuously differentiable since $s>1$. Thus by the mean-value theorem, there exist $\xi_1,\xi_2\in[0,\tau]$ such that
\begin{equation}\label{e141}
\begin{aligned}
\|I_4(&t_n)\|_{H^\frac12}=\|\textstyle\int_0^{\tfrac\tau2} f(\xi)d\xi-\tau f(\tfrac\tau2)+\int_{\tfrac\tau2}^\tau f(\xi)d\xi\|_{H^{\frac12}}
\leq\tau^2\|f^\prime(\xi_1)-f^\prime(\xi_2)\|_{H^{\frac12}}\\
&\lesssim\tau^2\|\fe^{-i(\tau/2-\xi_1)\langle\nabla\rangle}(\re u(t_n+\xi_1))^3)-\fe^{-i(\tau/2-\xi_2)\langle\nabla\rangle}(\re u(t_n+\xi_1))^3)\|_{H^\frac12}\\
&\quad+\tau^2\|\fe^{-i(\tau/2-\xi_2)\langle\nabla\rangle}(\re u(t_n+\xi_1))^3)-\fe^{-i(\tau/2-\xi_2)\langle\nabla\rangle}(\re u(t_n+\xi_2))^3)\|_{H^\frac12}\\
&\lesssim\tau^2\|(\fe^{-i(\xi_2-\xi_1)\langle\nabla\rangle}-1)\re (u(t_n+\xi_1))^3\|_{H^{\frac12}}\\
&\quad+\tau^2\|\re(u(t_n+\xi_1))^3-\re(u(t_n+\xi_2))^3\|_{H^{\frac12}},
\end{aligned}
\end{equation}
where the operator $\fe^{-i(\tau/2-\xi_2)\langle\nabla\rangle}$ has modulus 1 in Fourier spaces. By \eqref{cbe} and \eqref{boundu}, we have
\begin{equation}\label{e142}
\begin{aligned}
\|(\fe^{-i(\xi_2-\xi_1)\langle\nabla\rangle}-1)&(\re u(t_n+\xi_1))^3)\|_{H^{\frac12}}\lesssim\tau^{s-\frac12}\|u(t_n+\xi_1))^3\|_{H^s}\\
&\lesssim\tau^{s-\frac12}\|u(t_n+\xi_1))\|^3_{H^s}\lesssim\tau^{s-\frac12}.
\end{aligned}
\end{equation}

If we assume without loss of generality that $\xi_1\geq\xi_2$, then by \eqref{cbe}, \eqref{duh} and \eqref{boundu}, we get
\begin{equation}\label{e143}
\begin{aligned}
\|\re(u(t_n&+\xi_1))^3-\re(u(t_n+\xi_2))^3\|_{H^{\frac12}}\\
&\lesssim\|u(t_n+\xi_1)-u(t_n+\xi_2)\|_{H^\frac12}\sup\limits_{\xi\in[0,\tau]}\|u(t_n+\xi)\|^2_{H^s}\\
&\lesssim\|(\fe^{-i(\xi_1-\xi_2)\langle\nabla\rangle}-1)u(t_n+\xi_2)\|_{H^\frac12}+\tau\sup\limits_{\xi\in[0,\tau]}\|u(t_n+\xi)^3\|_{H^{-\frac12}}\\
&\lesssim\tau^{s-\frac12}\sup\limits_{\xi\in[0,\tau]}\|u(t_n+\xi)\|_{H^s}+\tau\sup\limits_{\xi\in[0,\tau]}\|u(t_n+\xi)\|^3_{H^s}\lesssim\tau^{s-\frac12}.
\end{aligned}
\end{equation}

The estimate \eqref{strangloc1} follows from combining \eqref{e13}, \eqref{e141}, \eqref{e142}, \eqref{e143}. This ends the proof.
\end{proof}

\subsection{Proof of Theorem \ref{theosob}}
First of all, by \eqref{kglie}, \eqref{kgstrang}, \eqref{clieloc} and \eqref{cstrangloc}, we can write the global errors of Lie and Strang respectively as
\begin{equation}\label{clieglo}
\begin{aligned}
p_n&:=u_n-u(t_n)=\Phi(u_{n-1})-\Phi(u(t_{n-1}))+\Phi(u(t_{n-1}))-u(t_n)\\
&=\fe^{-i\tau\langle\nabla\rangle}p_{n-1}-i\tau\langle\nabla\rangle^{-1}\fe^{-i\tau\langle\nabla\rangle}((\re u_{n-1})^3-(\re u(t_{n-1}))^3)+\mathcal{E}_t^1(u(t_{n-1})),
\end{aligned}
\end{equation}
and
\begin{equation}\label{cstranglo}
\begin{aligned}
q_n&:=v_n-u(t_n)=\Phi(v_{n-1})-\Psi(u(t_{n-1}))+\Psi(u(t_{n-1}))-u(t_n)\\
&=\fe^{-i\tau\langle\nabla\rangle}q_{n-1}-i\tau\langle\nabla\rangle^{-1}\fe^{-i\tau\langle\nabla\rangle/2}((\re \fe^{-i\tau\langle\nabla\rangle/2}u_{n-1})^3\\
&\quad-(\re\fe^{-i\tau\langle\nabla\rangle/2}u(t_{n-1}))^3)+\mathcal{E}_t^2(u(t_{n-1})).
\end{aligned}
\end{equation}
Then, we can prove Theorem~\ref{theosob} as follows.
\begin{proof}
For $s\in(\frac12,\frac56]$, by \eqref{multsob}, \eqref{lieloc1}, \eqref{boundu}, \eqref{clieglo}, we deduce that
\begin{align*}
\|p_n\|_{H^\frac12}&\leq\|p_{n-1}\|_{H^\frac12}+\tau\|(\re(u(t_{n-1})+p_{n-1}))^3-(\re u(t_{n-1}))^3\|_{H^{-\frac12}}+C_T\tau^{3s-\frac12}\\
&\leq(1+C_T\tau)\|p_{n-1}\|_{H^\frac12}+C_T\tau\|p_{n-1}\|^2_{H^\frac12}+\tau\|p_{n-1}\|^3_{H^\frac12}+C_T\tau^{3s-\frac12},
\end{align*}
where $C_T$ is a constant that depends on $s,~\|u_0\|_{H^s}$ and $T$. Since $p_0=0$, by the discrete Gronwall lemma together with a priori assumption and a bootstrap argument (see for example, \cite[pp.~2233]{WangZhao}), we thus conclude \eqref{liesob}.

For the Strang splitting, by \eqref{multsob}, \eqref{strangloc1}, \eqref{cstranglo} with $s\in(1,\frac32]$, we have
\begin{align*}
\|q_n\|_{H^\frac12}&\leq\tau\|\fe^{-3i\tau\langle\nabla\rangle/2}(\re(u(t_{n-1})+q_{n-1}))^3-\fe^{-3i\tau\langle\nabla\rangle/2}(\re u(t_{n-1}))^3\|_{H^{-\frac12}}\\
&\quad+\|q_{n-1}\|_{H^\frac12}+C_T\tau^{s+\frac32}\\
&\leq(1+C_T\tau)\|q_{n-1}\|_{H^\frac12}+C_T\tau\|q_{n-1}\|^2_{H^\frac12}+\tau\|q_{n-1}\|^3_{H^\frac12}+C_T\tau^{s+\frac32}.
\end{align*}
The estimate \eqref{strangsob} then follows similarly as $q_0=0$. This ends the proof.
\end{proof}

\begin{remark}\label{remsob}
For $s\in(\frac12,1]$, the Strang splitting method may exhibit less regular convergence behavior (see also Section~\ref{sectnumexp}). This is the reason why we did not include the case $s\in(\frac12,1]$ for the Strang splitting in Theorem~\ref{theosob}.
\end{remark}

%%%%%%%%%%%%%%%%%%%%%%
\section{Error estimates in Bourgain space}\label{sectbourg}
%%%%%%%%%%%%%%%%%%%%%%
The previous error analysis in the Sobolev space allows the lowest possible regularity of the solution to be $\frac12$, while the well-posedness only requires $s>\frac14$ \cite{Charxiv,Coweb}. To fill in the gap in between where the solution is truly very rough, we will work out the error estimates of the splitting schemes \eqref{kglie} and \eqref{kgstrang} utilizing the Bourgain space technique that will be specially developed for NKG in this section.

\subsection{A Bourgain framework}\label{2DBourg}
We first recall the definition and some important properties of the (discrete) Bourgain space \cite{Jiima,Ostjems}.

For a function $w(t,x)$ on $\mathbb{R}\times\mathbb{T}^2$, let $\tilde{w}(\sigma,k)$ stand for its time-space Fourier transform, i.e.,
$$
\tilde{w}(\sigma,k):=\int_{\mathbb{R}\times\mathbb{T}^2}w(t,x)\fe^{-i\sigma t-i\langle k,x\rangle}dx dt,
$$
where $\sigma\in\mathbb{R},~k\in\mathbb{Z}^d$, and $\langle\cdot,\cdot\rangle$ denotes the inner product in $\mathbb R^2$. 
The Bourgain space $X^{s, b}= X^{s,b}(\mathbb{R}\times\mathbb{T}^2)$ consists of functions with finite norm of the following form
\begin{equation}\label{contidef}
\|w\|_{X^{s,b}}:=\|\fe^{it\langle\nabla\rangle }w\|_{H_t^bH_x^s}=\|\langle k\rangle^s\langle\sigma+\langle k\rangle\rangle^b\tilde{w}(\sigma,k)\|_{ L^2l^2(\mathbb{R}\times\mathbb{Z}^2)}.
\end{equation}
To establish the multilinear estimate, i.e., the counterpart of \eqref{multsob} in discrete Bourgain spaces, we need to introduce the following projection operator/filter defined through the Fourier space, which is technically typical among the Bourgain settings \cite{Jisiam,Ostjems,Roufocm}:
\begin{equation}\label{filter}
\Pi_\tau:=\chi\left(\frac{-i\nabla}{\tau^{-1}}\right),
\end{equation}
where $\chi$ is the characteristic function of the square $[-1,1)^2$. Note that the filter depends on the size of $\tau$. 
With \eqref{filter}, we shall define a projected version of the equation \eqref{kgb}:
\begin{equation}\label{kgp}
u^\tau_t=-i\langle\nabla\rangle u^\tau-i\langle\nabla\rangle^{-1}\Pi_\tau(\re\Pi_\tau u^\tau)^3,\quad u^\tau(0,x)=\Pi_\tau u_0(x)=\Pi_\tau z_0(x)+i\langle\nabla\rangle^{-1}\Pi_\tau z_1(x);
\end{equation}
a projected version of the Lie splitting scheme:
\begin{equation}\label{plie}
u_{n+1}=\Phi^\tau(u_n)=\fe^{-i\tau\langle\nabla\rangle}u_n-i\tau\langle\nabla\rangle^{-1}\fe^{-i\tau\langle\nabla\rangle}\Pi_\tau(\re\Pi_\tau u_n)^3;
\end{equation}
and a projected version of the Strang splitting scheme:
\begin{equation}\label{pstrang}
u_{n+1}=\Psi^\tau(u_n)=\fe^{-i\tau\langle\nabla\rangle}u_n-i\tau\langle\nabla\rangle^{-1}\fe^{-i\tau\langle\nabla\rangle/2}\Pi_\tau(\re\Pi_\tau\fe^{-i\tau\langle\nabla\rangle/2}u_n)^3.
\end{equation} 
We then introduce some basic properties of $\Pi_\tau$, where the proof is rather straightforward from the definition \eqref{filter}: for any $s'\geq s,~b\in\mathbb{R}$ and $v$ defined in the corresponding spaces, we have
\begin{subequations}\label{projection prop}
\begin{align}
\label{pitau1}\|(I-\Pi_\tau)v\|_{X^{s,b}}+\|\Pi_\tau v\|_{X^{s,b}}&\lesssim\|v\|_{X^{s,b}},\\
\label{pitau2}\|(I-\Pi_\tau)v\|_{X^{s,b}}&\lesssim\tau^{s'-s}\|v\|_{X^{s',b}},\\
\label{pitau3}\|\Pi_\tau v\|_{X^{s',b}}&\lesssim\tau^{s-s'}\|v\|_{X^{s,b}}.
\end{align}
\end{subequations}

With the help of the Bourgain space \eqref{contidef}, we establish the following theorem in this section.
\begin{theorem}[Convergence in Bourgain]\label{theobourg}
For $s\in(\frac{11}{40},\frac{51}{40}]$ and the initial data $z_0(x)\in H^s(\mathbb{T}^2),$ $z_1(x)\in H^{s-1}(\mathbb{T}^2)$, let $(z,z_t)\in X^{s,b_0}\times X^{s-1,b_0}$ be the exact solution to equation \eqref{kg} on $[0,T]$ for $b_0\in(\frac12,\frac58)$, and denote $u=z+i\langle\nabla\rangle^{-1}z_t$ (i.e., $u$ is the solution to the equivalent equation \eqref{kgb}). Furthermore, we denote by $u_n,~v_n$ the sequences defined by the filtered Lie splitting method \eqref{plie} and the filtered Strang splitting method \eqref{pstrang}, respectively. Then, we have the following error estimate: there exist constants $\tau_0,~C_T>0$ such that for any time step size $\tau\in(0,\tau_0]$ and $0\leq n\tau\leq T$,
\begin{subequations}
\begin{align}
\label{liebourg}
\mbox{Lie: }\quad\|u(t_n)-u_n\|_{H^\frac{11}{40}}\leq C_T\tau^{s-\frac{11}{40}},\\
\label{strangbourg}
\mbox{Strang: }\quad\|u(t_n)-v_n\|_{H^\frac{11}{40}}\leq C_T\tau^{s-\frac{11}{40}},
\end{align}
\end{subequations}
where $C_T$ depends on $T$ and $\tau_0$ but is independent of $n$ and $\tau$.
\end{theorem}

Some well-known properties of Bourgain spaces are recalled as follows.

\begin{lemma}\label{contiprop}
For $\eta\in\mathcal{C}_c^\infty(\mathbb{R})$, and some space-time functions $F,~u_0,~u_1,~u_2,~u_3$ temporally supported in $[-2T,2T]$ for any $T>0$, we have that
\begin{subequations}\label{lm bourgain prop}
\begin{align}
\label{contiini}\|\eta(t)\fe^{-it\langle\nabla\rangle}u_0\|_{X^{s,b}}&\lesssim\|u_0\|_{H^s(\mathbb{T}^2)},\quad s\in\mathbb{R},~b\in\mathbb{R},~u_0\in H^s(\mathbb{T}^2);\\
\label{contib-1}\Big\|\int_0^t\fe^{-i(t-t')\langle\nabla\rangle}F(t')dt'\Big\|_{X^{s,b}}&\lesssim T^{1-b-b'}\|F\|_{X^{s,-b'}},\  s\in\mathbb{R},~b\in(\tfrac12,1],~b'\in(0,\tfrac12),~b+b'\leq1;\\
\label{contilinf}\|u_1\|_{L^\infty H^s}&\lesssim\|u_1\|_{X^{s,b}},\quad s\in\mathbb{R},~b>\tfrac12;\\
\label{contimult1}\|v_1v_2v_3\|_{X^{s-1,-b}}&\lesssim\|u_1\|_{X^{s,b}}\|u_2\|_{X^{s,b}}\|u_3\|_{X^{s,b}},\quad s>\tfrac{11}{40},~b>\tfrac38;\\
\label{contimult2}\|v_1v_2v_3\|_{X^{-\frac{29}{40},-b}}&\lesssim\|u_1\|_{X^{\frac{11}{40},b}}\|u_2\|_{X^{s,b}}\|u_3\|_{X^{s,b}},\quad s>\tfrac{11}{40},~b>\tfrac38;
\end{align}
\end{subequations}
where $v_j=u_j$ or $v_j=\overline{u_j}$, $j=1,2,3$.
\end{lemma}

We emphasize that these estimates are uniform for $T$, and we refer to \cite[Section~2.6]{Taoams} and \cite{Jiima} for the proofs of \eqref{contiini}--\eqref{contilinf}. The proof of the crucial nonlinear estimates~\eqref{contimult1} and \eqref{contimult2} is very similar to that of the discrete counterparts given in Section~\ref{sectprobourg}, which will be discussed later in Section~\ref{remhighd}.
We now start again with the Cauchy problem for \eqref{kgb}.

\begin{proposition}\label{theoexist2}
Let $s>\frac{11}{40}$ and $u_0\in H^s(\mathbb{T}^2)$. Then, there exist $T^*\in(0,+\infty]$ and a unique solution $u\in X^{s,b_0}\subset\mathcal{C}([0,T],H^s(\mathbb{T}^2))$ on $[0,T]$ for every $T\in(0,T^*)$ satisfying \eqref{kgb} and $b_0\in(\frac12,\frac58)$.
\end{proposition}
\begin{proof}
We consider the Duhamel form of \eqref{kgb}
$
v=F(v),
$
where 
\begin{equation}\label{duhunp} 
F(v)=\eta(t)\fe^{-it\langle\nabla\rangle}u_0-i\langle\nabla\rangle^{-1}\eta(t)\int_0^t\fe^{-i(t-\vartheta)\langle\nabla\rangle}(\re v(\vartheta))^3d\vartheta.
\end{equation}
Here $\eta\in[0,1]$ is supported in $[-2T_1,2T_1]$ and $\eta\equiv1$ on $[0,T_1]$ with $T_1\in(0,T]$ to be determined. By choosing $\frac38<b'<1-b_0<\frac12$ (so that $b'<\frac12<b_0$), denoting $\varepsilon_0=1-b_0-b'>0$ and utilizing \eqref{contiini}, \eqref{contib-1}, we have
$$
\|F(v)\|_{X^{s,b_0}}\lesssim\|u_0\|_{H^s}+T_1^{\varepsilon_0}\|v^3\|_{X^{s-1,-b'}},
$$
Then by using \eqref{contimult1}, we thus deduce that
$$
\|F(v)\|_{X^{s,b_0}}\leq C\|u_0\|_{H^s}+CT_1^{\varepsilon_0}\|v\|_{X^{s,b_0}}^3.
$$
By the same argument, we can get
$$
\|F(v_1)-F(v_2)\|_{X^{s,b_0}}\leq 3CT_1^{\varepsilon_0}R^2\|v_1-v_2\|_{X^{s,b_0}},
$$
for all $v_1,v_2$ satisfying $\|v_i\|_{X^{s,b_0}}\leq R$. Note that in the above estimates, $C$ is always independent of $T_1\in (0,T]$. The assertion then follows similarly as Proposition~\ref{theoexist}. Note that the resulting solution $u$ is defined globally but supported in $[-2T,2T]$ (see also in \cite{Jiima,Ostmcom}), and it is the solution to \eqref{kgb} only on $[0,T]$. 
\end{proof}

The next step is to study the solution of the projected equation \eqref{kgp}.

\begin{proposition}\label{theou-utau}
Let $u_0\in H^s$ for some $s>\frac{11}{40}$ and $u\in X^{s,b_0}$ the solution of \eqref{kgb} on $[0,T]$ given by Proposition~\ref{theoexist2}. Then, for $\tau$ sufficiently small, there exists a unique solution $u^\tau$ of \eqref{kgp} on $[0,T]$ such that $u^\tau\in X^{s,b_0}$. Moreover, we have
\begin{subequations}
\begin{align}
\label{diffutau}\|u-u^\tau\|_{X^{\frac{11}{40},b_0}}&\leq C_T\tau^{s-\frac{11}{40}},\\
\label{boundut}\|u^\tau\|_{X^{s,b_0}}&\leq C_T,
\end{align}
\end{subequations}
for some $C_T>0$ independent of $\tau$.
\end{proposition}
\begin{proof}
The existence and uniqueness of $u^\tau$ can be proven by a similar argument as given in Proposition~\ref{theoexist2}. We will prove $u^\tau$ is valid in $[0,T]$ with estimates \eqref{diffutau} and \eqref{boundut}. First of all, set
\begin{equation}\label{duhproj}
F^\tau(v)(t)=\eta(t)\fe^{-it\langle\nabla\rangle}\Pi_\tau u_0-i\langle\nabla\rangle^{-1}\eta(t)\int_0^t\fe^{-i(t-\vartheta)\langle\nabla\rangle}\Pi_\tau(\re\Pi_\tau v(\vartheta))^3d\vartheta,
\end{equation}
where $\eta\in[0,1]$ is supported in $[-2T_1,2T_1]$ and $\eta\equiv1$ on $[0,T_1]$ with $T_1\in(0,T]$ to be determined. Denote $u$ and $u^\tau$ the fixed-point of $F$ and $F^\tau$, respectively (the existence is guaranteed by Proposition~\ref{theoexist2}), where $F$ is defined in \eqref{duhunp}. Taking $b_1\in\big(\frac38,1-b_0)$, $\varepsilon_0=1-b_0-b_1>0,~s_1\in(\frac{11}{40},s)$, by \eqref{projection prop} and \eqref{lm bourgain prop}, noting that $u^\tau=\Pi_\tau u^\tau$, we have
\begin{align*}
\|u&-u^\tau\|_{X^{\frac{11}{40},b_0}}\lesssim\|(I-\Pi_\tau)u_0\|_{H^\frac{11}{40}}+T_1^{\varepsilon_0}\|(\re u)^3-\Pi_\tau(\re\Pi_\tau u^\tau)^3\|_{X^{-\frac{29}{40},-b_1}}\\
&\lesssim \tau^{s-\frac{11}{40}}\|u_0\|_{H^s}+T_1^{\varepsilon_0}\big(\|(I-\Pi_\tau)(\re u)^3\|_{X^{-\frac{29}{40},b_1}}+\|\Pi_\tau((\re u)^3-(\re\Pi_\tau u)^3)\|_{X^{-\frac{29}{40},b_1}}\\
&\quad+\|\Pi_\tau((\re\Pi_\tau u)^2\re\Pi_\tau(u-u^\tau))\|_{X^{-\frac{29}{40},b_1}}+\|\Pi_\tau((\re\Pi_\tau u-\re u^\tau)^2\re\Pi_\tau u)\|_{X^{-\frac{29}{40},b_1}}\\
&\quad+\|\Pi_\tau(\re\Pi_\tau u-\re u^\tau)^3\|_{X^{-\frac{29}{40},b_1}}\big)\\
&\lesssim \tau^{s-\frac{11}{40}}+T_1^{\varepsilon_0}(\tau^{s-\frac{11}{40}}\|u\|^3_{X^{s,b_1}}+\|u\|^2_{X^{s,b_1}}\|(I-\Pi_\tau)u\|_{X^{\frac{11}{40},b_1}}+\|u\|^2_{X^{s,b_1}}\|u-u^\tau\|_{X^{\frac{11}{40},b_1}}\\
&\quad+\|u\|_{X^{s,b_1}}\|u-u^\tau\|_{X^{\frac{11}{40},b_1}}\|\Pi_\tau(u-u^\tau)\|_{X^{s_1,b_1}}+\|u-u^\tau\|_{X^{\frac{11}{40},b_1}}\|\Pi_\tau(u-u^\tau)\|^2_{X^{s_1,b_1}}\big).
\end{align*}
Since $b_1<\frac12<b_0$, this yields
\begin{align*}
\|u-u^\tau\|_{X^{\frac{11}{40},b_0}}&\leq C_T\tau^{s-\frac{11}{40}}+C_TT_1^{\varepsilon_0}(\|u-u^\tau\|_{X^{\frac{11}{40},b_0}}+\tau^{\frac{11}{40}-s_1}\|u-u^\tau\|^2_{X^{\frac{11}{40},b_0}}\\
&\quad+\tau^{\frac{11}{20}-2s_1}\|u-u^\tau\|^3_{X^{\frac{11}{40},b_0}}).
\end{align*}
Noting that the initial data of $u-u^\tau$ is $(I-\Pi_\tau)u_0$, which satisfies
$$
\|(I-\Pi_\tau)u_0\|_{H^\frac{11}{40}}\lesssim\tau^{s-\frac{11}{40}},
$$
and that the choice of $C_T$ is independent of $T_1$, we thus conclude \eqref{diffutau} on $[0,T_1]$ by choosing $T_1$ sufficiently small such that $C_TT_1^{\varepsilon_0}\leq\frac12$ and $\|u-u^\tau\|_{X^{\frac{11}{40},b_0}}\lesssim\tau^{\frac{s+s_1}2-\frac{11}{40}}$. Then, the estimate follows from a similar iteration as given in \cite[Lemma 3.5]{Jimcom}.

It remains to prove \eqref{boundut}. In fact, by \eqref{pitau3}, we get
$$
\|u^\tau\|_{X^{s,b_0}}\leq\|\Pi_\tau(u-u^\tau)\|_{X^{s,b_0}}+\|\Pi_\tau u\|_{X^{s,b_0}}\leq\tau^{\frac{11}{40}-s}\|u-u^\tau\|_{X^{\frac{11}{40},b_0}}+\|u\|_{X^{s,b_0}}\leq C_T.
$$
This ends the proof.
\end{proof}

We now give the definition as well as some basic properties of the discrete (in time) Bourgain space. For more details, see also \cite{Jiima,Ostjems}.

We first take $\{u_n(x)\}_{n\in\mathbb{Z}}$ to be a sequence of functions on the torus $\mathbb{T}^2$, with its ``time-space" Fourier transform:
$$
\widetilde{u_n}(\sigma,k)=\tau\sum\limits_{m\in\mathbb{Z}}\widehat{u_m}(k)\fe^{-im\tau\sigma},\qquad
\mbox{where }\ 
\widehat{u_m}(k)=\int_{\mathbb{T}^2}u_m(x)\fe^{-i\langle k,x\rangle}dx.
$$
In this framework, $\widetilde{u_{n}}$ is a $2\pi/\tau$ periodic function in $\sigma$. Then, the discrete Bourgain space $X_\tau^{s,b}$ can be defined with the norm
\begin{equation}\label{discdef}
\|u_n\|_{X^{s,b}_\tau}=\|\fe^{in\tau\langle\nabla\rangle}u_n\|_{H_\tau^bH_x^s}=\|\langle k\rangle^s\langle d_\tau(\sigma+\langle k\rangle)\rangle^b\widetilde{u_n}(\sigma,k)\|_{ L^2l^2(\mathbb{T}\times\mathbb{Z}^2)},
\end{equation}
where $d_\tau(\sigma)=\frac{\fe^{i\tau\sigma}-1}{\tau}$. Note that we always have $|d_\tau|\lesssim\tau^{-1}$, and for any fixed $u_n$, the norm is an increasing function for both $s$ and $b$.

We directly get from the definition \eqref{discdef} the elementary properties that for any $s\geq s^\prime$ and $b\geq b^\prime$, 
\begin{subequations}
\begin{align}
\label{shift}\sup\limits_{\delta\in[-4,4]}\|\fe^{i\tau\delta\langle\nabla\rangle} u_n\|_{X^{s,b}_\tau}&\lesssim\|u_n\|_{X_\tau^{s,b}},\\
\label{b-bp}\|u_n\|_{X^{s,b}_\tau}&\lesssim\tau^{b^\prime-b}\|u_n\|_{X_\tau^{s,b^\prime}},\\
\label{s-sp}\|\Pi_\tau u_n\|_{X^{s,b}_\tau}&\lesssim\tau^{s^\prime-s}\|u_n\|_{X_\tau^{s^\prime,b}}.
\end{align}
\end{subequations}

Next, we will give the counterparts of Lemma~\ref{contiprop}.

\begin{lemma}\label{discprop}
For $s\in\mathbb{R},~\eta\in\mathcal{C}_c^\infty(\mathbb{R}),~\tau\in(0,1]$ and $u_n$ supported on $-2T\leq n\tau\leq2T$, we have that
\begin{align}
\label{discini}\|\eta(n\tau)\fe^{-in\tau\langle\nabla\rangle}u_0\|_{X_\tau^{s,b}}&\lesssim_{\eta,b}\|u_0\|_{H^s(\mathbb{T}^2)},\quad b\in\mathbb{R},~u_0\in H^s,\\
\label{disclinf}\|u_n\|_{l_\tau^\infty H^s}&\lesssim_b\|u_n\|_{X^{s,b}_\tau},\quad b>\tfrac12,\\
\label{discb-1}\|\tau\sum\limits_{m=0}^n\fe^{-i(n-m)\tau\langle\nabla\rangle}u_m(x)\|_{X^{s,b}_\tau}&\lesssim T^{1-b-b^\prime}\|u_n\|_{X^{s,-b'}_\tau},\quad b\in(\tfrac12,1),~b'\in(0,\tfrac12),~b+b'\leq1.
\end{align}
\end{lemma}

Note that all these estimates are uniform in $\tau$. 
The proof of this lemma is largely the same as the estimates given in \cite[Section~3]{Ostjems}. Therefore, we omit the details for simplicity.

We close this subsection by introducing a vital tool stated as the following theorem, which is the discrete counterpart of~\eqref{contimult1} and~\eqref{contimult2}. It plays a crucial role for the coming error estimates of the nonlinear part in NKG, and can be applied to other related models in the future.

\begin{theorem}[Discrete Bourgain estimate for wave-type]\label{theomultbourg}
For any $s>\frac{11}{40},~s'>\frac38,~b>\frac38$, we have
\begin{subequations}
\begin{align}
\label{l4strich}\|\Pi_\tau u_n\|_{l^4_\tau L^4}&\lesssim\|u_n\|_{X_\tau^{s',b}},\\
\label{multbourg1}\|\Pi_\tau(\Pi_\tau v_{1,n}\Pi_\tau v_{2,n}\Pi_\tau v_{3,n})\|_{X_\tau^{s-1,-b}}&\lesssim\|u_{1,n}\|_{X_\tau^{s,b}}\|u_{2,n}\|_{X_\tau^{s,b}}\|u_{3,n}\|_{X_\tau^{s,b}},\\
\label{multbourg2}\|\Pi_\tau(\Pi_\tau v_{1,n}\Pi_\tau v_{2,n}\Pi_\tau v_{3,n})\|_{X_\tau^{-\frac{29}{40},-b}}&\lesssim\|u_{1,n}\|_{X_\tau^{\frac{11}{40},b}}\|u_{2,n}\|_{X_\tau^{s,b}}\|u_{3,n}\|_{X_\tau^{s,b}},
\end{align}
\end{subequations}
where $v_{j,n}=u_{j,n}$ or $v_{j,n}=\overline{u_{j,n}}$, $j=1,2,3$, and $u_n$, $u_{1,n}$, $u_{2,n}$, $u_{3,n}$ are functions in the corresponding spaces.
\end{theorem}

We postpone the proof of this important theorem to Section~\ref{sectprobourg}.

\subsection{Local error estimates}

In this part, we estimate the local error of the filtered splitting schemes \eqref{plie} and \eqref{pstrang}. Here, we shall only estimate the error of the filtered Lie splitting method. The convergence behavior and analysis for the Strang case are similar. Note that, in the Bourgain setting, due to the filter $\Pi_\tau$, it is not possible to achieve better convergence results using the Strang splitting than Lie in the low-regularity situation \cite{Ostmcom}. For the case $s>\frac{51}{40}$, we recommend Sobolev settings without filters, as these provide a better convergence order.

To estimate the error in discrete Bourgain space, we need to prove the boundedness of $u^\tau$ in the target space. We the following lemma:

\begin{lemma}
For any $s>\frac{11}{40}$, $b_0\in(\frac12,\frac58)$, if $u_0\in H^s(\mathbb{T}^2)$, we have
\begin{equation}\label{discboundu}
\sup_{\vartheta\in[-4\tau,4\tau]}\|u^\tau(t_n+\vartheta)\|_{X_\tau^{s,b_0}}\leq C_T,
\end{equation}
where $u^\tau$ is the exact solution to \eqref{kgp}.
\end{lemma}
\begin{proof}
By setting $f=\fe^{-it\langle\nabla\rangle}v$ and $f_n(x)=f(n\tau,x)$, it suffices to prove that
$$
\|f_n\|_{H_\tau^bL_x^2(\mathbb{Z}\times\mathbb{T}^2)}\lesssim\|f\|_{H_t^bL_x^2},
$$
then the lemma follows from a combination of \eqref{boundut}. This estimate can be found, for instance, in \cite[Lemma 3.4]{Roupaa}. Note that since we allow $b>\frac12$, we do not need $\tau$ to compensate the gain of $b$.
\end{proof}

Next, similar to \eqref{clieloc}, we shall calculate the local error:

\begin{equation}\label{blieloc}
\begin{aligned}
\mathcal{E}_t^3(u^\tau(t_n))&=\Phi^\tau(u^\tau(t_n))-u^\tau(t_{n+1})\\
&=i\langle\nabla\rangle^{-1}\fe^{-i\tau\langle\nabla\rangle}\int_0^\tau(\fe^{i\xi\langle\nabla\rangle}-1)\Pi_\tau(\re\Pi_\tau  u^\tau(t_n+\xi))^3d\xi\\
&\quad+i\langle\nabla\rangle^{-1}\fe^{-i\tau\langle\nabla\rangle}\int_0^\tau\left[ \Pi_\tau(\re\Pi_\tau u^\tau(t_n+\xi))^3-\Pi_\tau(\re\Pi_\tau u^\tau(t_n))^3\right]d\xi\\
&=:J_1(t_n)+J_2(t_n).
\end{aligned}
\end{equation}

We now give the local error estimate result.

\begin{proposition}
For $s\in(\frac{11}{40},\frac{51}{40}]$, $b_0\in(\frac12,\frac58)$, we denote by $\mathcal{E}_t^3(u^\tau(t_n))$ the local error of \eqref{plie} as we calculated in \eqref{blieloc}, where $u^\tau$ is the exact solution of \eqref{kgp} as in Proposition~\ref{theou-utau}. Then, we have
\begin{equation}\label{lieloc2}
\|\mathcal{E}_t^3(u^\tau(t_n))\|_{X_\tau^{\frac{11}{40},b_0-1}}\leq C_T\tau^{s+\frac{29}{40}}.
\end{equation}
\end{proposition}
\begin{proof}
By \eqref{shift}, \eqref{multbourg1} and \eqref{discboundu}, we have
\begin{equation}\label{e31}
\begin{aligned}
\|J_1(t_n)\|_{X_\tau^{\frac{11}{40},b_0-1}}&\lesssim\tau\sup\limits_{\xi\in[0,\tau]}\|(\fe^{i\xi\langle\nabla\rangle}-1)\Pi_\tau(\re\Pi_\tau u^\tau(t_n+\xi))^3\|_{X_\tau^{-\frac{29}{40},b_0-1}}\\
&\lesssim\tau\sup\limits_{\xi\in[0,\tau]}\|(\xi\langle\nabla\rangle)^{s-\frac{11}{40}}\Pi_\tau(\re\Pi_\tau u^\tau(t_n+\xi))^3\|_{X^{-\frac{29}{40},b_0-1}_\tau}\\
&\lesssim\tau^{s+\frac{29}{40}}\sup\limits_{\xi\in[0,\tau]}\|\Pi_\tau(\re\Pi_\tau u^\tau(t_n+\xi))^3\|_{X^{s-1,b_0-1}_\tau}\\
&\lesssim\tau^{s+\frac{29}{40}}\sup\limits_{\xi\in[0,\tau]}\|u^\tau(t_n+\xi)\|^3_{X^{s,b_0}_\tau}\lesssim\tau^{s+\frac{29}{40}}.
\end{aligned}
\end{equation}

Similar to \cite{Jisiam,Jiima, Ostjems}, by interpolating \eqref{disclinf} and the trivial embedding $l_\tau^2H^s=X_\tau^{s,0}$, we have for any $s\in\mathbb{R},~b>\frac14$ and $v_n$ that
\begin{equation}\label{discl4}
\|v_n\|_{l_\tau^4H^s}\lesssim\|v_n\|_{X_\tau^{s,b}},
\end{equation}
and a dual version
\begin{equation}\label{discl34}
\|v_n\|_{X_\tau^{s,-b}}\lesssim\|v_n\|_{l_\tau^{\frac43}H^s}.
\end{equation}

Then, noting $b_0<\frac34$, by \eqref{duhproj}, \eqref{b-bp}, \eqref{s-sp}, \eqref{multbourg2}, and \eqref{discboundu}, we get for $s\leq\frac{31}{80}$ that
\begin{equation}\label{e32}
\begin{aligned}
\|J_2&(t_n)\|_{X_\tau^{\frac{11}{40},b_0-1}}\leq\tau\sup\limits_{\xi\in[0,\tau]}\|\Pi_\tau(\Pi_\tau u^\tau(t_n+\xi))^3-\Pi_\tau(\Pi_\tau u^\tau(t_n))^3\|_{X_\tau^{-\frac{29}{40},b_0-1}}\\
&\lesssim\tau\sup\limits_{\xi\in[0,\tau]}\|u^\tau(t_n+\xi)-u^\tau(t_n)\|_{X_\tau^{\frac{11}{40},\frac38+\varepsilon}}\sup\limits_{\xi\in[0,\tau]}\|\Pi_\tau u^\tau(t_n+\xi)\|^2_{X_\tau^{s,b_0}}\\
&\lesssim\tau\!\sup\limits_{\xi\in[0,\tau]}\|(\fe^{-i\xi\langle\nabla\rangle}-1)u^\tau(t_n)\|_{X_\tau^{\frac{11}{40},b_0}}+\tau^{\frac54-2\varepsilon}\!\sup\limits_{\xi\in[0,\tau]}\|\Pi_\tau(\Pi_\tau u^\tau(t_n+\xi))^3\|_{X_\tau^{-\frac{29}{40},-\frac38-\varepsilon}}\\
&\lesssim\tau^{s+\frac{29}{40}}\|u^\tau(t_n)\|_{X_\tau^{s,b_0}}+\tau^{\frac54-2\varepsilon}\sup\limits_{\xi\in[0,\tau]}\|\Pi_\tau u^\tau(t_n+\xi)\|^3_{X_\tau^{s,b_0}}\lesssim\tau^{s+\frac{29}{40}},
\end{aligned}
\end{equation}
where $\varepsilon\in(0,\frac1{16})$ can be taken arbitrarily small. For the case $s>\frac{31}{80}$, by using \eqref{duhproj}, \eqref{l4strich}, \eqref{discl4}, \eqref{discl34}, H\"older's inequality and Sobolev embedding theorem, we obtain
\begin{equation}\label{e33}
\begin{aligned}
\|J_2&(t_n)\|_{X_\tau^{\frac{11}{40},b_0-1}}\leq\tau\sup\limits_{\xi\in[0,\tau]}\|\Pi_\tau(\Pi_\tau u^\tau(t_n+\xi))^3-\Pi_\tau(\Pi_\tau u^\tau(t_n))^3\|_{l_\tau^\frac43H^{-\frac{29}{40}}}\\
&\lesssim\tau\sup\limits_{\xi\in[0,\tau]}\|u^\tau(t_n+\xi)-u^\tau(t_n)\|_{l_\tau^4L^\frac{80}{29}}\sup\limits_{\xi\in[0,\tau]}\|\Pi_\tau u^\tau(t_n+\xi)\|^2_{l_\tau^4L^4}\\
&\lesssim\tau\sup\limits_{\xi\in[0,\tau]}\|u^\tau(t_n+\xi)-u^\tau(t_n)\|_{l_\tau^4 H^\frac{11}{40}}\sup\limits_{\xi\in[0,\tau]}\|\Pi_\tau u^\tau(t_n+\xi)\|^2_{X_\tau^{s,b_0}}\\
&\lesssim\tau\sup\limits_{\xi\in[0,\tau]}\|(\fe^{-i\xi\langle\nabla\rangle}-1)u^\tau(t_n)\|_{X_\tau^{\frac{11}{40},b_0}}+\tau^2\sup\limits_{\xi\in[0,\tau]}\|(\Pi_\tau u^\tau(t_n+\xi))^3\|_{l_\tau^4 H^{-\frac{29}{40}}}\\
&\lesssim\tau^{s+\frac{29}{40}}\|u^\tau(t_n)\|_{X_\tau^{s,b_0}}+\tau^2\sup\limits_{\xi\in[0,\tau]}\|\Pi_\tau u^\tau(t_n+\xi)\|_{l_\tau^4L^4}\|\Pi_\tau u^\tau(t_n+\xi)\|^2_{l_\tau^\infty L^\frac{160}{49}}\\
&\lesssim\tau^{s+\frac{29}{40}}+\tau^2\sup\limits_{\xi\in[0,\tau]}\|\Pi_\tau u^\tau(t_n+\xi)\|^3_{X_\tau^{s,b_0}}\lesssim\tau^{s+\frac{29}{40}}.
\end{aligned}
\end{equation}
We conclude the proposition by combining \eqref{e31}, \eqref{e32} and \eqref{e33}.
\end{proof}

\subsection{Proof of Theorem~\ref{theobourg}}

In this part, we will estimate the global error of the filtered Lie splitting method \eqref{plie}, and prove Theorem~\ref{theobourg}. We only consider the filtered Lie splitting (i.e., to prove \eqref{liebourg}) as an example, while \eqref{strangbourg} can be proved similarly. First of all, we shall write the global error as follows:
\begin{equation}\label{blieglo}
\begin{aligned}
e_n&=u_n-u^\tau(t_n)\\
&=\fe^{-i\tau\langle\nabla\rangle}\big(u_{n-1}-u^\tau(t_{n-1})\big)-i\tau\langle\nabla\rangle^{-1}\fe^{-i\tau\langle\nabla\rangle}\big(\Omega(u_{n-1})-\Omega(u^\tau(t_{n-1}))\big)\\
&\quad+\mathcal{E}_t^3(u^\tau(t_{n-1}))\\
&=-i\tau\langle\nabla\rangle^{-1}\sum\limits_{k=1}^{n-1}\fe^{-i(n-k)\tau\langle\nabla\rangle}\big(\Omega(u_k)-\Omega(u^\tau(t_k))\big)+\sum\limits_{k=0}^{n-1}\fe^{-i(n-k-1)\tau\langle\nabla\rangle}\mathcal{E}_t^3(u^\tau(t_k)),
\end{aligned}
\end{equation}
where
$
\Omega(u)=\Pi_\tau(\re\Pi_\tau u)^3.
$ 
Now, we are ready to prove Theorem~\ref{theobourg}.
\begin{proof}
By \eqref{discb-1}, \eqref{lieloc2}, \eqref{blieglo}, we have for $e_n$ defined as \eqref{blieglo} on $[0,T_1]$ and supported in $[-2T_1,2T_1]$ that
$$
\|e_n\|_{X_\tau^{\frac{11}{40},b_0}}\leq C_TT_1^{\varepsilon_0}\|\Pi_\tau\big(\!\re\Pi_\tau(e_n+u^\tau(t_n))\big)^3-\Pi_\tau\big(\!\re\Pi_\tau u^\tau(t_n)\big)^3\|_{X_\tau^{-\frac{29}{40},-b_1}}+C_T\tau^{s-\frac{11}{40}},
$$
where $b_0\in(\frac12,\frac58),~b_1\in(\frac38,1-b_0),~\varepsilon_0=1-b_0-b_1$ and $T_1$ to be determined. Thus by \eqref{s-sp}, \eqref{multbourg2}, we obtain
\begin{align*}
\|e_n\|_{X_\tau^{\frac{11}{40},b_0}}&\leq C_TT_1^{\varepsilon_0}\big(\|u^\tau(t_n)\|^2_{X_\tau^{s,b_0}}\|e_n\|_{X_\tau^{\frac{11}{40},b_0}}+\|u^\tau(t_n)\|_{X_\tau^{s,b_0}}\|e_n\|_{X_\tau^{\frac{11}{40},b_0}}\|\Pi_\tau e_n\|_{X_\tau^{s_1,b_0}}\\
&\quad+\|\Pi_\tau e_n\|^2_{X_\tau^{s_1,b_0}}\|e_n\|_{X_\tau^{\frac{11}{40},b_0}}\big)+C_T\tau^{s-\frac{11}{40}}.
\end{align*}
This yields that for $s_1\in(\frac{11}{40},s)$,
$$
\|e_n\|_{X_\tau^{\frac{11}{40},b_0}}\leq C_TT_1^{\varepsilon_0}\big(\|e_n\|_{X_\tau^{\frac{11}{40},b_0}}+\tau^{\frac{11}{40}-s_1}\|e_n\|_{X_\tau^{\frac{11}{40},b_0}}^2+\tau^{\frac{11}{20}-2s_1}\|e_n\|_{X_\tau^{\frac{11}{40},b_0}}^3\big)+C_T\tau^{s-\frac{11}{40}},
$$
since $\Pi_\tau e_n=e_n$. Noting that the choice of $C_T$ depends on $s,~b_1,~\|u_0\|_{H^s}$ but is independent of $T_1$, so we can take $T_1$ sufficiently small such that $C_TT_1^{\varepsilon_0}\leq\frac12$ to get
$$
\|e_n\|_{X_\tau^{\frac{11}{40},b_0}}\leq C_T\tau^{s-\frac{11}{40}}
$$
on $[0,T_1]$. Reiterating this estimate on $[T_1,2T_1]$, $[2T_1,3T_1]$ and so on (see also~\cite{Jiima}), we arrive at
$$
\|e_n\|_{X_\tau^{\frac{11}{40},b_0}}\leq C_T\tau^{s-\frac{11}{40}}
$$
on $[0,T]$.

Finally, by \eqref{diffutau}, \eqref{disclinf}, we obtain
\begin{align*}
\|u_n-u(t_n)\|_{H^\frac{11}{40}}&\leq\|u_n-u^\tau(t_n)\|_{l_\tau^\infty H^\frac{11}{40}}+\|u-u^\tau\|_{L^\infty H^\frac{11}{40}}\\
&\leq C'_T\|u_n-u^\tau(t_n)\|_{X_\tau^{\frac{11}{40},b_0}}+C'_T\|u-u^\tau\|_{X^{\frac{11}{40},b_0}}\leq C_T\tau^{s-\frac{11}{40}},
\end{align*}
which concludes \eqref{liebourg}.  This ends the proof.
\end{proof}

%%%%%%%%%%%%%%%%%%%%%%%%%%%%%%%
\section{Nonlinear estimates and generalizations}\label{sectproof}
%%%%%%%%%%%%%%%%%%%%%%%%%%%%%%%

In this section, we will prove the multilinear estimates given in Lemma~\ref{theomultsob} and Theorem~\ref{theomultbourg}. We will reformulate them and prove a more general version so that the estimates work not only for $d=2$ but also for other dimensions. Some more general discussions will be made in the end.

\subsection{Proof of Lemma \ref{theomultsob}}\label{sectprosob}

In this part, we will reformulate Lemma \ref{theomultsob} as follows and then prove it.

\begin{proposition}
For any function $u_i\in H^{s_i}(\mathbb{T}^d),~i=1,2,3$, suppose $r\in(-\tfrac d2,\tfrac d2)$, $s_i\in(r,\tfrac d2),~s_1+s_2+s_3\geq r+d$. Then, we have the following estimate:
\begin{equation}\label{multsobnd}
\|u_1u_2u_3\|_{H^r}\lesssim\|u_1\|_{H^{s_1}}\|u_2\|_{H^{s_2}}\|u_3\|_{H^{s_3}}.
\end{equation}
\end{proposition}
\begin{proof}
Without loss of generality, we assume $s_1\leq s_2\leq s_3$, and $s_1+s_2+s_3=r+d$ (otherwise we can take the smaller $s_i$). Since $r+d>\tfrac d2>s_3$, we must have $s_2>0$. Then, we split the discussion into the following three cases:
\begin{enumerate}
\item $s_1>r\geq0$.

By the Kato--Ponce inequality \cite{Bourgdie,Katocpam}, we have for $p_i\in(2,+\infty),~1\leq i\leq9$ that
\begin{equation}\label{kp1}
\left\{\begin{aligned}
\|u_1u_2u_3\|_{H^r}&\lesssim\|u_1\|_{W^{r,p_1}}\|u_2\|_{L^{p_2}}\|u_3\|_{L^{p_3}}+\|u_1\|_{L^{p_4}}\|u_2\|_{W^{r,p_5}}\|u_3\|_{L^{p_6}}\\
&\quad+\|u_1\|_{L^{p_7}}\|u_2\|_{L^{p_8}}\|u_3\|_{W^{r,p_9}}\\
\tfrac1{p_1}+\tfrac1{p_2}+\tfrac1{p_3}&=\tfrac12,\quad \tfrac1{p_4}+\tfrac1{p_5}+\tfrac1{p_6}=\tfrac12,\quad \tfrac1{p_7}+\tfrac1{p_8}+\tfrac1{p_9}=\tfrac12.
\end{aligned}\right.
\end{equation}
Note that with $j=1,4,7$,
$$
\frac{s_1}d+\frac{s_2}d+\frac{s_3}d-\frac rd=1=\frac12-\frac1{p_j}+\frac12-\frac1{p_{j+1}}+\frac12-\frac1{p_{j+2}},
$$
where we have
$$
\frac{s_i-r}d,\ \frac{s_i}d\in(0,\frac12),\quad i=1,2,3.
$$
Thus, if we take
$$
\left\{\begin{aligned}
p_1=\frac1{\frac12-\frac{s_1-r}d},\quad p_2=\frac1{\frac12-\frac{s_2}d},\quad p_3=\frac1{\frac12-\frac{s_3}d},\\
p_4=\frac1{\frac12-\frac{s_1}d},\quad p_5=\frac1{\frac12-\frac{s_2-r}d},\quad p_6=\frac1{\frac12-\frac{s_3}d},\\
p_7=\frac1{\frac12-\frac{s_1}d},\quad p_8=\frac1{\frac12-\frac{s_2}d},\quad p_9=\frac1{\frac12-\frac{s_3-r}d},
\end{aligned}\right.
$$
then by \eqref{kp1} and the Sobolev embedding theorem, we obtain \eqref{multsobnd}.

\item $s_1\geq0>r$.

By the dual version of the Sobolev embedding theorem, it suffices to prove
\begin{equation}\label{holder0}
\|u_1u_2u_3\|_{L^{p_0}}\lesssim\|u_1\|_{H^{s_1}}\|u_2\|_{H^{s_2}}\|u_3\|_{H^{s_3}},
\end{equation}
where $\frac1{p_0}-\frac12=-\frac rd$. By H\"older's inequality, we get
\begin{equation}\label{holder1}
\|u_1u_2u_3\|_{L^{p_0}}\lesssim\|u_1\|_{L^{p_1}}\|u_2\|_{L^{p_2}}\|u_3\|_{L^{p_3}},
\end{equation}
where $\tfrac1{p_1}+\tfrac1{p_2}+\tfrac1{p_3}=\tfrac1{p_0}$. This yields
$$
\frac12-\frac1{p_1}+\frac12-\frac1{p_2}+\frac12-\frac1{p_3}=1+\frac rd=\frac{s_1}d+\frac{s_2}d+\frac{s_3}d.
$$
We therefore set
$$
p_1=\frac1{\frac12-\frac{s_1}d},\quad p_2=\frac1{\frac12-\frac{s_2}d},\quad p_3=\frac1{\frac12-\frac{s_3}d},
$$
then \eqref{holder0} follows from \eqref{holder1} and the Sobolev embedding theorem.
\item $0>s_1>r$.

Take $u_0$ to be any function in $H^{-r}$. Then, \eqref{multsobnd} is equivalent to
$$
|\textstyle\int_{\mathbb{T}^d}u_0u_1u_2u_3dx|\lesssim\|u_0\|_{H^{-r}}\|u_1\|_{H^{s_1}}\|u_2\|_{H^{s_2}}\|u_3\|_{H^{s_3}}.
$$
Moreover, this is also equivalent to
$$
\|u_0u_2u_3\|_{H^{-s_1}}\lesssim\|u_0\|_{H^{-r}}\|u_2\|_{H^{s_2}}\|u_3\|_{H^{s_3}}.
$$
Note that $s_1+s_2=r+d-s_3>-\frac d2+d-\frac d2=0$, i.e., $s_2>-s_1$, we therefore have $-s_1\in(-\frac d2,\frac d2)$, $0<-s_1<-r,s_2,s_3<\frac d2$ and $-r+s_2+s_3=-s_1+d$. Then this case is the same as the first one.
\end{enumerate}
We conclude \eqref{multsobnd} by collecting the three cases.
\end{proof}

\begin{remark}\label{rem1d}
In fact, the Kato-Ponce inequality also allows us to have $p_i=\infty$ (except $i=1,5,9$). However, the difference is that the embedding $L^\infty\supset H^{\frac d2}$ is not true. Instead, we can only prove $L^\infty\supset H^{\frac d2+\varepsilon}$. From this, we can prove for example in the one-dimensional case that
$$
\|u_1u_2u_3\|_{H^r}\lesssim\|u_1\|_{H^{s_1}}\|u_2\|_{H^{s_2}}\|u_3\|_{H^{s_3}},\quad r<-\tfrac12,~s_1+s_2+s_3=\tfrac12,~s_i\in(-\tfrac12,\tfrac12).
$$
In particular, we can take $r=-\frac56$ and $s_1=s_2=s_3=\frac16$ to get the local well-posedness of \eqref{kgb} for $s=\frac16$ in the one-dimensional case.

In addition, we also have to note that at least two of the three variables $s_i$ need to be non-negative.
\end{remark}

\subsection{Proof of Theorem \ref{theomultbourg}}\label{sectprobourg}

Next, we prove Theorem~\ref{theomultbourg}. We shall focus on the two-dimensional case, and discuss the generalization to other dimensions in Section~\ref{remhighd}. 
We will first recall the Littlewood--Paley decomposition, then introduce several technical lemmas, and finally prove the desired estimates.

Let us first recall the Littlewood--Paley decomposition. We define for $\sigma\in I_{\tau}=[-\frac\pi\tau,\frac\pi\tau)$ and $m\in\mathbb{N}_0$ that
$P_m(\sigma)=\mathbbm{1}_{2^m\leq\langle\sigma\rangle<2^{m+1}\cap I_\tau}$.
$P_m$ also denotes the $\frac{2\pi}\tau$ periodic extension of this function. 
Similarly, we define localizations in spatial frequencies. We set
$Q_m(k)=\mathbbm{1}_{2^m\leq\langle k\rangle<2^{m+1}}$.
Moreover, we set $N=(N_0,N_1,N_2,N_3),~L=(L_0,L_1,L_2,L_3)$, where $N_j$ and $L_j$ are dyadic of the form
$N_j=2^{n_j},L_j=2^{l_j},~j=0,1,2,3$. We split $u_{j,n}$ as
$$
u_{j,n}=\sum\limits_{l_j,n_j\in\mathbb{N}}u_{j,n}^{L_jN_j}, \quad\mbox{with}\quad \widetilde{u_{j,n}^{L_jN_j}}(\sigma,k)=P_{l_j}(\sigma+\langle k\rangle)Q_{n_j}(k)\widetilde{u_{j,n}}(\sigma,k),\quad j=0,1,2,3.
$$
Then, we define
\begin{equation}\label{defsln}
S(L,N)=\Big|\tau\sum_n\int_{\mathbb{T}^2}\Pi_\tau v_{0,n}^{L_0N_0}\Pi_\tau v_{1,n}^{L_1N_1}\Pi_\tau v_{2,n}^{L_2N_2}\Pi_\tau v_{3,n}^{L_3N_3}dx\Big|,
\end{equation}
so that
\begin{equation}\label{ssln}
S=\Big|\tau\sum_n\int_{\mathbb{T}^2}\Pi_\tau v_{0,n}\Pi_\tau v_{1,n}\Pi_\tau v_{2,n}\Pi_\tau v_{3,n}dx\Big|\leq\sum_{L,N}S(L,N),
\end{equation}
where, similar to Theorem~\ref{theomultbourg}, $v_{j,n}=u_{j,n}$ or $v_{j,n}=\overline{u_{j,n}}$, $j=0,1,2,3$.

We then give a technical lemma.

\begin{lemma}\label{lemlocl4}
For any $p\in\mathbb{Z}^2$ with $|p|=M\lesssim\tau^{-1}$, $v_n$ satisfying $\widetilde{v_n}(\sigma,k)$ localized at $k-p\in[-N,N)^2$ and $\sigma+\langle k\rangle=\mathcal{O}(L)$, we have
\begin{subequations}
\begin{align}
\label{locl4v1}\|\Pi_\tau v_n\|_{l_\tau^4L^4}&\lesssim (M+N)^\frac18N^\frac14L^\frac38\|v_n\|_{l_\tau^2L^2},\\
\label{locl4v2}\|\Pi_\tau v_n\|_{l_\tau^4L^4}&\lesssim N^\frac12L^\frac14\|v_n\|_{l_\tau^2L^2}.
\end{align}
\end{subequations}
\end{lemma}
\begin{proof}
We first prove \eqref{locl4v2}. Without loss of generality, we assume $\|v_n\|_{l_\tau^2L^2}=1$. By Plancherel's theorem, it suffices to show that
$$
\Bigg\|\sum\limits_{\substack{k_{11}+k_{21}=k_1\\k_{12}+k_{22}=k_2}}\int_{\sigma_1+\sigma_2=\sigma}\widetilde{\Pi_\tau v_n}(\sigma_1,k_{11},k_{12})\widetilde{\Pi_\tau v_n}(\sigma_2,k_{21},k_{22})d\sigma_1\Bigg\|_{L^2l^2}\lesssim NL^\frac12.
$$
Note that we have by assumption $\sigma_i+\sqrt{k_{i1}^2+k_{i2}^2}\in E_L$, $i=1,2$ with $E_L\subset\bigcup_{m\in\mathbb{Z}}[\frac{2m\pi}{\tau}-L,\frac{2m\pi}{\tau}+L]$. Since $\sigma_1,~\sigma\in I_\tau$ and $k_{ij}$ satisfies $|k_{ij}|\leq\tau^{-1}$ due to $\Pi_\tau$, we actually have $E_L=\bigcup_{|m|\leq B}[\frac{2m\pi}{\tau}-L,\frac{2m\pi}{\tau}+L]$ with $B=\mathcal{O}(1)$, i.e., $|E_L|=\mathcal{O}(L)$. 

Following similarly to the Schr\"odinger case given in \cite[Lemma 8.1]{Jiima}, it suffices to show that
$$
|\Omega_0|=\big|\big\{(k_{11},k_{12}):\sigma+\textstyle\sqrt{k_{11}^2+k_{12}^2}+\textstyle\sqrt{(k_1-k_{11})^2+(k_2-k_{12})^2}\in E_{2L}\big\}\big|=\mathcal{O}(N^2),
$$
where all the numbers except $(k_{11},k_{12})$ are fixed. This is obvious since both $k_{11}$ and $k_{12}$ have at most $N$ choices. This concludes \eqref{locl4v2}.

For \eqref{locl4v1}, by a similar argument, it suffices to show that
$$
|\Omega_0|=\mathcal{O}(N\sqrt{(M+N)L}).
$$
If we denote $O$ as the origin of a plane and denote $F=(k_1,k_2)$, then our aim is to count the number of integer points $A=(k_{11},k_{12})$ such that the length of the piecewise linear segment $|OA|+|AF|$ is $\sigma_0+\mathcal{O}(L)$, where $\sigma_0=\sigma+\frac{2m\pi}\tau$ with proper $m$. Note that $\sigma_0$ may not be unique. Nevertheless, we can substitute every possible $\sigma_0$ in the following discussion to obtain the final result.

Since both $k_{11}$ and $k_{12}$ have only $\mathcal{O}(N)$ choices, we assume without loss of generality that $L<N$. We then split the problem into three cases:
\begin{enumerate}
\item $F$ is different from $O$, and $M\leq 5N$.

This means that $A$ lies on an ellipse with foci $O$ and $F$, and the length of major axis is $2a=\sigma_0+\mathcal{O}(L)$. If we denote by $a_1$ the largest possible semi-major axis while $a_2$ the smallest possible semi-major axis, then $a_1-a_2=\mathcal{O}(L)$. We also denote by $c$ the semi-focal length, i.e., $c=\frac12\sqrt{k_1^2+k_2^2}$. If we denote by $b_1$ and $b_2$ the largest and smallest possible semi-minor axis, then $b_1^2=a_1^2-c^2,~b_2^2=a_2^2-c^2$.

We know that $A$ lies in an elliptic ring, and note that the width of the most narrow part of this ring is $a_1-a_2=\mathcal{O}(L)\gtrsim1$, we thus get that
$$
|\Omega_0|\lesssim\pi a_1b_1-\pi a_2b_2\lesssim b_1(a_1-a_2)+a_2(b_1-b_2).
$$
Note that $a_1=\mathcal{O}(N)$, $a_2=\mathcal{O}(N)$, and
$$
(b_1-b_2)^2\leq b_1^2-b_2^2=a_1^2-a_2^2=(a_1+a_2)(a_1-a_2)=\mathcal{O}(NL),
$$
we thus have
$$
|\Omega_0|\lesssim b_1(a_1-a_2)+a_2(b_1-b_2)=\mathcal{O}(NL+N\sqrt{NL})=\mathcal{O}(\sqrt{N^3L}),
$$
which concludes the first case.

\item $F$ is different from $O$, and $M>5N$.

Here we use the same notation as in the first case. Here $A$ lies at the intersection of the elliptic ring as introduced in the first case, and a square with side length $N$. This is an area of two long bands, where each of them has a length of $\mathcal{O}(N)$ and width at most
$$
b_1-b_2\leq\sqrt{b_1^2-b_2^2}=\sqrt{a_1^2-a_2^2}=\sqrt{(a_1+a_2)(a_1-a_2)}=\mathcal{O}(\sqrt{ML}).
$$
Consequently, we also have
$
|\Omega_0|=\mathcal{O}(N\sqrt{ML})
$
for this case.

\item $F$ coincides with $O$.

This means that both $A=(k_{11},k_{12})$ and $A'=(-k_{11},-k_{12})$ satisfy $k-a\in[-N,N)^2$, which gives $|p|=\mathcal{O}(N)$. Then, the integer point $A$ lies in a ring with radius $O(N)$ and width $\mathcal{O}(L)$. The number of such points is $\mathcal{O}(N^2-(N-L)^2)=\mathcal{O}(NL)\leq\mathcal{O}(\sqrt{N^3L})$.
\end{enumerate}

We finally conclude \eqref{locl4v1} by collecting the three cases.
\end{proof}

Next, in order to reduce the regularity constraint given by \eqref{l4strich}, we will introduce the following lemma.

\begin{lemma}
For any $j_1,~j_2\in\{0,1,2,3\}$, we have
\begin{subequations}
\begin{align}
\label{j1j2}\big\|\Pi_\tau u_{j_1,n}^{L_{j_1}N_{j_1}}\Pi_\tau u_{j_2,n}^{L_{j_2}N_{j_2}}\big\|_{l_\tau^2L^2}&\lesssim N_{j_1}^\frac58N_{j_2}^\frac18\big\|u_{j_1,n}^{L_{j_1}N_{j_1}}\big\|_{X_\tau^{0,\frac38}}\big\|u_{j_2,n}^{L_{j_2}N_{j_2}}\big\|_{X_\tau^{0,\frac38}},\\
\label{j3j4}\big\|\Pi_\tau u_{j_1,n}^{L_{j_1}N_{j_1}}\Pi_\tau u_{j_2,n}^{L_{j_2}N_{j_2}}\big\|_{l_\tau^2L^2}&\lesssim N_{j_1}^\frac78\big\|u_{j_1,n}^{L_{j_1}N_{j_1}}\big\|_{X_\tau^{0,\frac38}}\big\|u_{j_2,n}^{L_{j_2}N_{j_2}}\big\|_{X_\tau^{0,\frac38}}.
\end{align}
\end{subequations}
\end{lemma}
\begin{proof}
We adapt the proof of NLS given in \cite[Lemma 8.3]{Jiima}. Without loss of generality, we assume $j_1=1,~j_2=2$ and $N_1\leq N_2$. For simplicity of notation, we denote $v_n=u_{1,n}^{L_1N_1},~w_n=u_{2,n}^{L_2N_2}$.

We first divide $\tilde{w}_n$ into $\mathcal{O}\big((\frac{N_2}{N_1})^2\big)$ localized functions, where the frequency of each is localized in a square of side length $N_1$. If we denote by $R_p$ the localization operator, where $p\in\mathbb{R}^2$ is the center of the square, we have by definition that
\begin{equation}\label{defrp}
\sum\limits_p\|R_p w_n\|^2_{l_\tau^2L^2}=\|w_n\|^2_{l_\tau^2L^2}.
\end{equation}

Moreover, since $k_1+k_2=p+\mathcal{O}(N_1)$, the sum $\sum\limits_p\Pi_\tau v_n\Pi_\tau R_p w_n$ is actually quasi-orthogonal. Thus by \eqref{locl4v1}, \eqref{defrp} and Cauchy--Schwarz's inequality, we deduce that
\begin{align*}
\|\Pi_\tau v_n&\Pi_\tau w_n\|_{l_\tau^2L^2}\lesssim\Big(\sum\limits_p\|\Pi_\tau v_n\Pi_\tau R_p w_n\|^2_{l_\tau^2L^2}\Big)^\frac12\lesssim\Big(\sum\limits_p\|\Pi_\tau v_n\|^2_{l_\tau^4L^4}\|\Pi_\tau R_p w_n\|^2_{l_\tau^4L^4}\Big)^\frac12\\
&\lesssim N_1^\frac58N_2^\frac18L_1^\frac38L_2^\frac38\Big(\|v_n\|^2_{l_\tau^2L^2}\Big)^\frac12\Big(\sum\limits_p\|R_p w_n\|^2_{l_\tau^2L^2}\Big)^\frac12\lesssim N_1^\frac58N_2^\frac18\|v_n\|_{X_\tau^{0,\frac38}}\|w_n\|_{X_\tau^{0,\frac38}},
\end{align*}
which concludes \eqref{j1j2}. Similarly, by \eqref{locl4v1}, \eqref{locl4v2}, \eqref{defrp} and Cauchy--Schwarz's inequality, we obtain
$$
\|\Pi_\tau v_n\Pi_\tau w_n\|_{l_\tau^2L^2}\lesssim N_1^\frac78L_1^\frac38L_2^\frac14\Big(\|v_n\|^2_{l_\tau^2L^2}\Big)^\frac12\Big(\sum\limits_p\|R_p w_n\|^2_{l_\tau^2L^2}\Big)^\frac12\lesssim N_1^\frac78\|v_n\|_{X_\tau^{0,\frac38}}\|w_n\|_{X_\tau^{0,\frac38}},
$$
which concludes \eqref{j3j4}. This ends the proof.
\end{proof}

We now move on to the proof of Theorem~\ref{theomultbourg}. 

\begin{proof}
We first prove \eqref{l4strich}. Again, we write
$$
u_n=\sum\limits_{l,m\in\mathbb{N}}u_n^{LN}, \qquad \widetilde{u_n^{LN}}(\sigma,k)=P_l(\sigma+\langle k\rangle)Q_m(k)\widetilde{u_n}(\sigma,k),
$$
where $L=2^l,~N=2^m$. Thus by \eqref{locl4v1} and Cauchy--Schwarz's inequality, we have
\begin{equation}\label{sumup}
\begin{aligned}
\|\Pi_\tau&u_n\|^2_{l_\tau^4L^4}\lesssim\Big(\sum\limits_{L,N}\|\Pi_\tau u_n^{LN}\|_{l_\tau^4L^4}\Big)^2\lesssim\Big(\sum\limits_{L,N}N^\frac38L^\frac38\|\Pi_\tau u_n^{LN}\|_{l_\tau^2L^2}\Big)^2\\
&\lesssim\Big(\sum\limits_{L,N}N^{\frac34-2s}L^{\frac34-2b}\Big)\Big(\sum\limits_{L,N}N^{2s}L^{2b}\|\Pi_\tau u_n^{LN}\|^2_{l_\tau^2L^2}\Big)\\&\sim\sum\limits_{L,N}\|\Pi_\tau u_n^{LN}\|^2_{X_\tau^{s,b}}\sim\|\Pi_\tau u_n\|^2_{X_\tau^{s,b}},
\end{aligned}
\end{equation}
which concludes \eqref{l4strich}.

For \eqref{multbourg1}, we assume without loss of generality that $s<\frac12$. The case that $s\geq\frac12$ can be proved by combining \eqref{multsob}, \eqref{discl4} and \eqref{discl34}. Take $u_{0,n}\in X_\tau^{1-s,b}$, and set $v_{0,n}=u_{0,n}$ or $v_{0,n}=\overline{u_{0,n}}$, then by \eqref{defsln} and \eqref{ssln}, it suffices to show
\begin{equation}\label{multeqv}
S\leq\sum\limits_{L,N}S(L,N)\lesssim\|u_{0,n}\|_{X_\tau^{1-s,b}}\|u_{1,n}\|_{X_\tau^{s,b}}\|u_{2,n}\|_{X_\tau^{s,b}}\|u_{3,n}\|_{X_\tau^{s,b}}.
\end{equation}
Without loss of generality, we assume $N_1\geq N_2\geq N_3$. Note that $N_0\lesssim N_1+N_2+N_3\lesssim N_1$ (otherwise $S(L,N)$ vanishes), then by \eqref{j1j2} and Cauchy--Schwarz's inequality, we have
\begin{equation}\label{slnest1}
\begin{aligned}
S(L,N)&\lesssim\|v_{0,n}^{L_0N_0}v_{1,n}^{L_1N_1}\|_{l_\tau^2L^2}\|v_{2,n}^{L_2N_2}v_{3,n}^{L_3N_3}\|_{l_\tau^2L^2}=\|u_{0,n}^{L_0N_0}u_{1,n}^{L_1N_1}\|_{l_\tau^2L^2}\|u_{2,n}^{L_2N_2}u_{3,n}^{L_3N_3}\|_{l_\tau^2L^2}\\
&\lesssim\|u_{0,n}^{L_0N_0}\|_{X_\tau^{\frac58,\frac38}}\|u_{1,n}^{L_1N_1}\|_{X_\tau^{\frac18,\frac38}}\|u_{2,n}^{L_2N_2}\|_{X_\tau^{\frac18,\frac38}}\|u_{3,n}^{L_3N_3}\|_{X_\tau^{\frac58,\frac38}}.
\end{aligned}
\end{equation}
Similarly, by \eqref{j1j2}, \eqref{j3j4} and Cauchy--Schwarz's inequality, we get
\begin{equation}\label{slnest2}
S(L,N)\lesssim\|u_{0,n}^{L_0N_0}\|_{X_\tau^{\frac78,\frac38}}\|u_{1,n}^{L_1N_1}\|_{X_\tau^{0,\frac38}}\|u_{2,n}^{L_2N_2}\|_{X_\tau^{\frac18,\frac38}}\|u_{3,n}^{L_3N_3}\|_{X_\tau^{\frac58,\frac38}}.
\end{equation}
If we denote $s_1=\frac13(s-\frac{11}{40})$, note that $\frac15-2s_1=\frac{23}{60}-\frac23s>0$, then by interpolating \eqref{slnest1} with strength $\frac35$ and \eqref{slnest2} with strength $\frac25$, and using Cauchy--Schwarz's inequality, we obtain
\begin{align*}
S&\lesssim\sum\limits_{L,N}N_0^\frac{29}{40}N_1^\frac3{40}N_2^\frac18N_3^\frac58\|u_{0,n}^{L_0N_0}\|_{X_\tau^{0,\frac38}}\|u_{1,n}^{L_1N_1}\|_{X_\tau^{0,\frac38}}\|u_{2,n}^{L_2N_2}\|_{X_\tau^{0,\frac38}}\|u_{3,n}^{L_3N_3}\|_{X_\tau^{0,\frac38}}\\
&\lesssim\sum\limits_{L,N}N_0^{\frac{29}{40}-4s_1}N_1^{\frac3{40}+4s_1+\frac15-2s_1}N_2^{\frac{11}{40}+s_1}N_3^{\frac{11}{40}+s_1}\\
&\quad\times\|u_{0,n}^{L_0N_0}\|_{X_\tau^{0,\frac38}}\|u_{1,n}^{L_1N_1}\|_{X_\tau^{0,\frac38}}\|u_{2,n}^{L_2N_2}\|_{X_\tau^{0,\frac38}}\|u_{3,n}^{L_3N_3}\|_{X_\tau^{0,\frac38}}\\
&\lesssim\big(\sum\limits_{L,N}(N_0N_1)^{-2s_1}(N_2N_3)^{-4s_1}(L_0L_1L_2L_3)^{\frac34-2b}\big)^\frac12\\
&\quad\times\Big(\sum\limits_{L,N}\|u_{0,n}^{L_0N_0}\|^2_{X_\tau^{1-s,b}}\|u_{1,n}^{L_1N_1}\|^2_{X_\tau^{s,b}}\|u_{2,n}^{L_2N_2}\|^2_{X_\tau^{s,b}}\|u_{3,n}^{L_3N_3}\|^2_{X_\tau^{s,b}}\Big)^\frac12\\
&\lesssim\|u_{0,n}\|_{X_\tau^{1-s,b}}\|u_{1,n}\|_{X_\tau^{s,b}}\|u_{2,n}\|_{X_\tau^{s,b}}\|u_{3,n}\|_{X_\tau^{s,b}},
\end{align*}
which concludes \eqref{multbourg1}.

As for \eqref{multbourg2}, since it is not symmetric, we assume $\sigma\in S_3$ such that $N_{\sigma(1)}\geq N_{\sigma(2)}\geq N_{\sigma(3)}$. Note that we still have $N_0\lesssim N_{\sigma(1)}$. Thus similarly, we find
\begin{align*}
S&\lesssim\sum\limits_{L,N}N_0^\frac{29}{40}N_{\sigma(1)}^\frac3{40}N_{\sigma(2)}^\frac18N_{\sigma(3)}^\frac58\|u_{0,n}^{L_0N_0}\|_{X_\tau^{0,\frac38}}\|u_{1,n}^{L_1N_1}\|_{X_\tau^{0,\frac38}}\|u_{2,n}^{L_2N_2}\|_{X_\tau^{0,\frac38}}\|u_{3,n}^{L_3N_3}\|_{X_\tau^{0,\frac38}}\\
&\lesssim\sum\limits_{L,N}N_0^{\frac{29}{40}-s_1}N_{\sigma(1)}^{\frac3{40}+s_1+\frac15-2s_1}N_{\sigma(2)}^{\frac{11}{40}+s_1}N_{\sigma(3)}^{\frac{11}{40}+s_1}
\\&\quad\times\|u_{0,n}^{L_0N_0}\|_{X_\tau^{0,\frac38}}\|u_{1,n}^{L_1N_1}\|_{X_\tau^{0,\frac38}}\|u_{2,n}^{L_2N_2}\|_{X_\tau^{0,\frac38}}\|u_{3,n}^{L_3N_3}\|_{X_\tau^{0,\frac38}}\\
&\lesssim\sum\limits_{L,N}N_0^{\frac{29}{40}-s_1}N_1^{\frac{11}{40}-s_1}N_2^{\frac{11}{40}+s_1}N_3^{\frac{11}{40}+s_1}\|u_{0,n}^{L_0N_0}\|_{X_\tau^{0,\frac38}}\|u_{1,n}^{L_1N_1}\|_{X_\tau^{0,\frac38}}\|u_{2,n}^{L_2N_2}\|_{X_\tau^{0,\frac38}}\|u_{3,n}^{L_3N_3}\|_{X_\tau^{0,\frac38}}\\
&\lesssim\|u_{0,n}\|_{X_\tau^{\frac{29}{40},b}}\|u_{1,n}\|_{X_\tau^{\frac{11}{40},b}}\|u_{2,n}\|_{X_\tau^{s,b}}\|u_{3,n}\|_{X_\tau^{s,b}}.
\end{align*}
This ends the proof.
\end{proof}

\subsection{Generalization and further discussions}\label{remhighd}

Following closely the lines of the proof given in Section~\ref{sectprobourg}, one can prove the continuous nonlinear estimates \eqref{contimult1} and \eqref{contimult2}. Note that, in the continuous case, since $\sigma\in\mathbb{R}$ and the functions are no longer periodic in $\sigma$, cancellations occur only at $0$ (unlike in the discrete case, where cancellations occur at $2k\pi\tau^{-1}$ with $k\in\mathbb{Z}$ \cite{Roufocm}), so the filter $\Pi_\tau$ is not necessary.

The constant $11/40$ given in the section above is not optimal, where in principle one can prove similar estimates for $s>\frac14$ \cite{Coweb,Lindjfa}. Unlike in the Schr\"odinger case, the estimate \eqref{locl4v1} given in Lemma~\ref{lemlocl4} is not independent of $M$, which introduces additional difficulties. By introducing a more involved Strichartz framework, the authors in \cite{Lindjfa} overcame this issue and proved the local well-posedness for $s>\frac14$. These techniques far exceed the scope of this paper.

The analysis of this paper is restricted to the parameters $m=1,~\lambda=-1$. This can be easily extended to the focusing NKG, i.e., $m=1,~\lambda=1$. Moreover, following the lines of the proof given in this paper, one can even extend to nonlinear wave equations, where $m=0$ and $\langle\nabla\rangle$ is not invertible. Note that, as introduced in \cite{Baomcom}, the schemes \eqref{kglie} and \eqref{kgstrang} are actually the same as the methods generated by the following splitting technique:
$$
\left\{\begin{aligned}
    &u_t=v,\\&v_t=\Delta u-mu;
\end{aligned}\right.
\qquad\qquad \rm{and} \qquad\qquad
\left\{\begin{aligned}
    &u_t=0,\\&v_t=\lambda u^3.
\end{aligned}\right.
$$
So there is no singularity from the method point of view. In order to simplify the discussions in Bourgain spaces, we omit the discussions on the case $m=0$ in this paper. However, since the value of $m$ does not make differences to the local well-posedness theories \cite{Coweb}, by analyzing the zero frequency term of $z_t$ in addition, the case $m=0$ follows similarly.

For $d=1$, the Strichartz estimate \eqref{l4strich} in principle (i.e., for integrable functions) works for $s>0$ with $b<\frac12$ sufficiently close to $\frac12$. However, this estimate cannot provide integrability of the nonlinearity due to the lack of Strichartz smoothing, so multilinear estimates \eqref{contimult1}, \eqref{contimult2}, \eqref{multbourg1}, \eqref{multbourg2} do not work for $s>0$~\cite{Charxiv,Coweb}. Instead, the constant $\frac16$ given by the Sobolev estimate \eqref{multsobnd} is critical for the one-dimensional case (see also Remark~\ref{rem1d}).

For $d\geq3$, the $l_\tau^4L^4$ Strichartz estimate is different. Specifically, we can prove
\begin{equation}\label{l4nd}
\|\Pi_\tau u_n\|_{l^4_\tau L^4}\lesssim\|u_n\|_{X_\tau^{s,b}},
\end{equation}
for $s>\frac{d-1}4$ and $b<\frac12$ sufficiently close to $\frac12$. The proof is quite similar to Lemma~\ref{lemlocl4}, the difference is that we first choose $k_{i1},~3\leq i\leq d$ arbitrarily without changing the square sum $\sum\limits_{i=3}^dk_{i1}^2$, then the problem reduces to the three-dimensional case, i.e., counting the volume of an ellipsoidal ring. If we use the same notations as in Lemma~\ref{lemlocl4}, by similar arguments as given in the lemma, we get
$$
|\Omega_0|\lesssim a_1b_1^2-a_2b_2^2=b_1^2(a_1-a_2)+a_2(b_1^2-b_2^2)=b_1^2(a_1-a_2)+a_2(a_1+a_2)(a_1-a_2)=\mathcal{O}(N^2L),
$$
and the case where $M>5N$ follows similarly. This yields
\begin{equation}\label{interp1}
\|\Pi_\tau v_n\|_{l_\tau^4L^4}\lesssim (M+N)^\frac14 N^\frac{d-2}4L^\frac12\|v_n\|_{l_\tau^2L^2}.
\end{equation}
Moreover, by a similar argument as given in Lemma~\ref{lemlocl4}, we have
\begin{equation}\label{interp2}
\|\Pi_\tau v_n\|_{l_\tau^4L^4}\lesssim N^\frac{d}4L^\frac14\|v_n\|_{l_\tau^2L^2}.
\end{equation}
Substituting $u_n^{LN}$ into these estimates, interpolating \eqref{interp1} and \eqref{interp2}, and synthesizing the results via \eqref{sumup}, we finally arrive at \eqref{l4nd} noting also that $M=0$ for $u_n^{LN}$.

For the nonlinear estimates, by \eqref{interp1}, \eqref{interp2} and using similar techniques given in this subsection, we can prove
\begin{align*}
\|\Pi_\tau(\Pi_\tau v_{1,n}\Pi_\tau v_{2,n}\Pi_\tau v_{3,n})\|_{X_\tau^{s-1,-b}}&\lesssim\|u_{1,n}\|_{X_\tau^{s,b}}\|u_{2,n}\|_{X_\tau^{s,b}}\|u_{3,n}\|_{X_\tau^{s,b}},\\
\|\Pi_\tau(\Pi_\tau v_{1,n}\Pi_\tau v_{2,n}\Pi_\tau v_{3,n})\|_{X_\tau^{-\frac12,-b}}&\lesssim\|u_{1,n}\|_{X_\tau^{s,b}}\|u_{2,n}\|_{X_\tau^{s,b}}\|u_{3,n}\|_{X_\tau^{\frac12,b}},
\end{align*}
where $s>\frac d2-1$. Unlike the two-dimensional case, the constant $\frac d2-1$ is sharp for local well-posedness \cite{Coweb}. Due to these inequalities, we can estimate the error in $H^\frac12$. Moreover, using similar techniques introduced in Section~\ref{sectbourg}, we can prove that the convergence order of the method is $s-\frac12$ in $H^\frac12$.

%%%%%%%%%%%%%%%%%%%%%%
\section{Numerical Experiments}\label{sectnumexp}
%%%%%%%%%%%%%%%%%%%%%%

In this section, we numerically illustrate the convergence results given in Theorem~\ref{theosob} and Theorem~\ref{theobourg}. We shall generate initial data with different regularities to test the error of the numerical schemes under the proved norms. To be specific, we take the initial data
$$
z_0(x)+\langle\nabla\rangle^{-1}z_1(x)=\sum_{k_1,k_2=-N/2}^{N/2-1}\langle k\rangle^{-s-\frac d2-\varepsilon}\tilde{f}_k\fe^{ikx}\in H^s,
$$
where $k=(k_1,k_2)$, $\varepsilon>0$ can be taken arbitrarily small, and $N$ is the number of spatial grid points. Here we numerically take $\varepsilon=0$. Moreover, $\tilde{f}_k,~\tilde{g}_k$ are all random variables that are uniformly distributed in the interval $[0,1]$. We employ the standard Fourier pseudospectral method (FFT) for spatial discretization and fix $T=1$ as the final time.

\begin{figure}
\begin{center}
\subfigure{\includegraphics[width=0.49\textwidth]{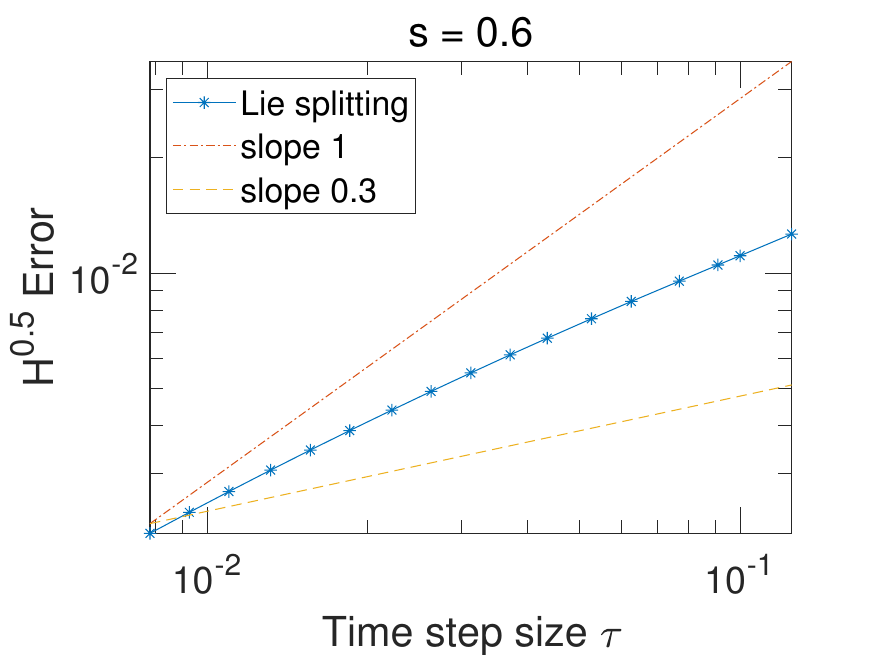}}
\subfigure{\includegraphics[width=0.49\textwidth]{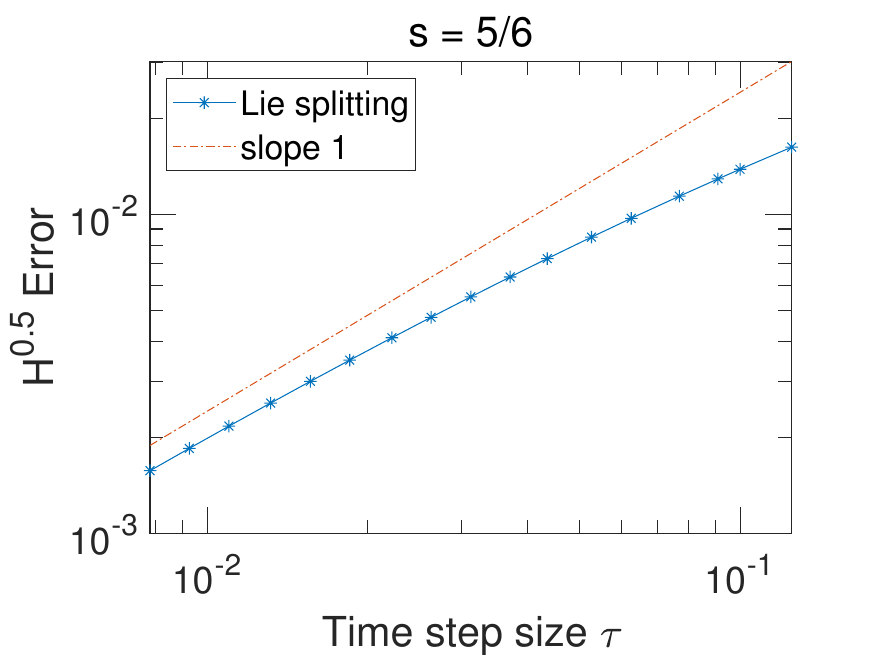}}
\end{center}
\caption{$H^{\frac12}$ temporal error of the Lie splitting scheme for the two-dimensional NKG with rough initial data $z_0+\langle\nabla\rangle^{-1}z_1\in H^s$.  
Left: $s=0.6$; \ right: $s=5/6$.\label{fig1}}
\end{figure}

In Figure~\ref{fig1}, we performed numerical experiments for the Lie splitting scheme \eqref{kglie}. Here we normalized the initial data $\|z_0+\langle\nabla\rangle^{-1}z_1\|_{H^\frac12}=1$, and took Fourier modes $N=(2^{13},2^{13})$. For the reference solution, we used a small time step size $\tau=2^{-12}$. 
From the figure on the left, we can see that the numerical results are better than predicted in Theorem~\ref{theosob}. Nevertheless, we can clearly observe an order reduction. Here, possible reasons for the superconvergence include the finite Fourier series having a smoothing effect or some implicit regularization from global accumulation. From the figure on the right, we can see that the numerical method converges with the order $1$, which fits our theorem very well.

\begin{figure}
\begin{center}
\subfigure{\includegraphics[width=0.49\textwidth]{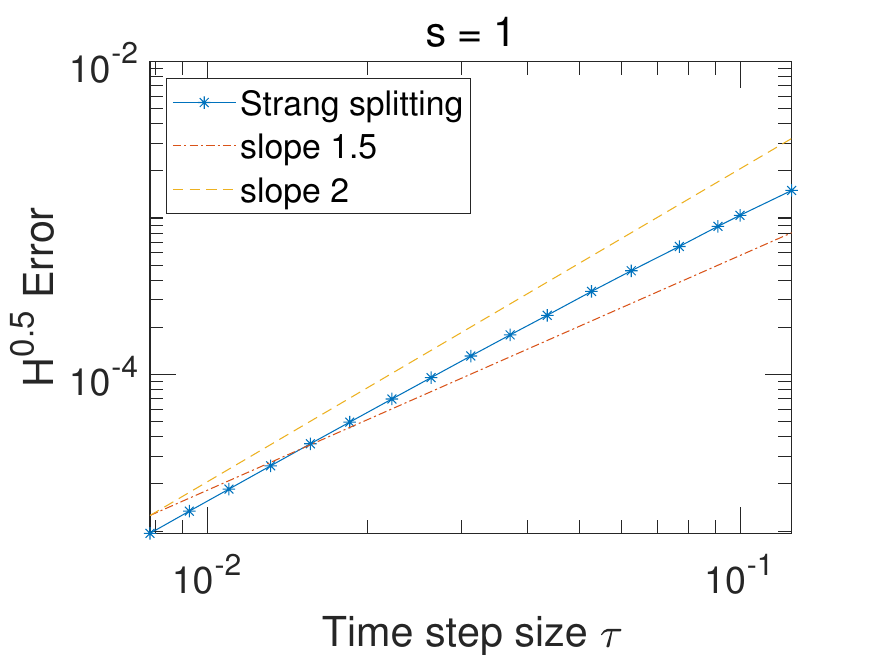}}
\subfigure{\includegraphics[width=0.49\textwidth]{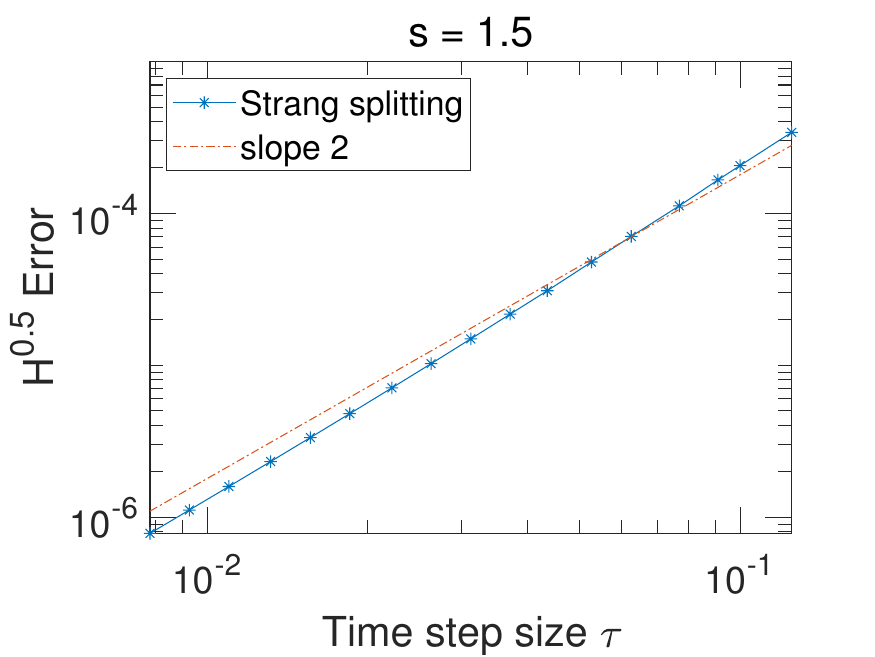}}
\end{center}
\caption{$H^{\frac12}$ temporal error of the Strang splitting scheme for the two-dimensional NKG with rough initial data $z_0+\langle\nabla\rangle^{-1}z_1\in H^s$.  
Left: $s=1$; \ right: $s=1.5$.\label{fig2}}
\end{figure}

In Figure~\ref{fig2}, we performed numerical experiments for the Strang splitting scheme \eqref{kgstrang}. Here we again normalized the initial data $\|z_0+\langle\nabla\rangle^{-1}z_1\|_{H^\frac12}=1$, and took Fourier modes $N=(2^{13},2^{13})$. For the reference solution, we used a small time step size $\tau=2^{-12}$. 
From the figure on the left, we can see that the error converges with order a bit more than $s+\frac12$, while an order reduction also exists, which is similar to Figure~\ref{fig1}. From the figure on the right, we can see that the numerical method converges with order $2$, which fits our theorem very well.

\begin{figure}
\begin{center}
\subfigure{\includegraphics[width=0.49\textwidth]{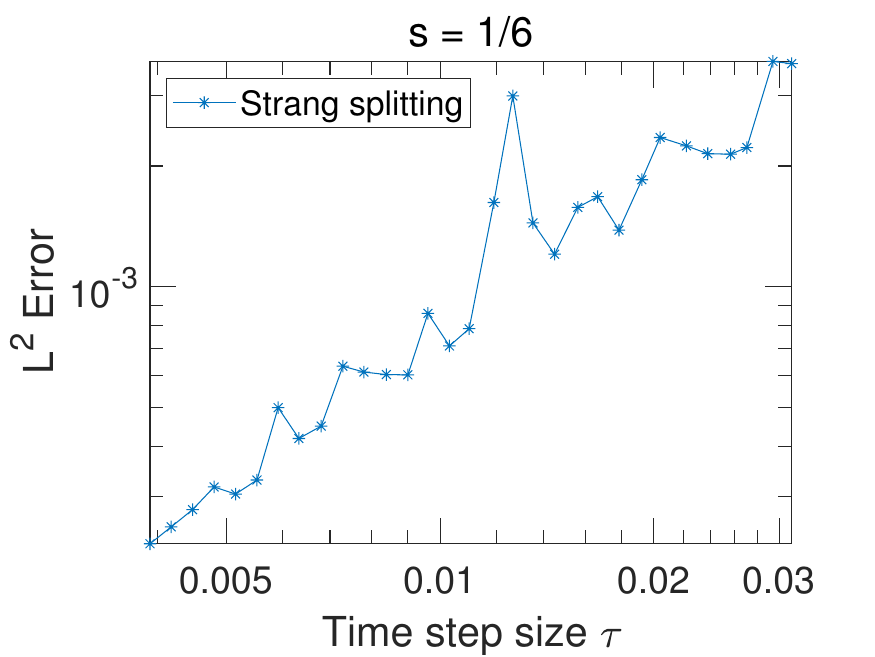}}
\subfigure{\includegraphics[width=0.49\textwidth]{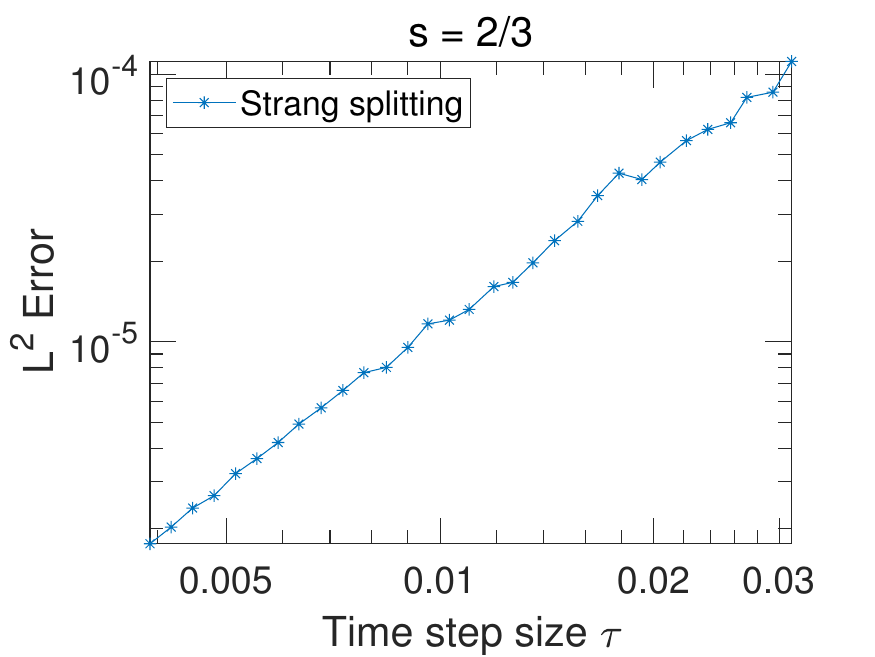}}
\end{center}
\caption{$L^2$ temporal error of the Strang splitting scheme for the one-dimensional NKG with rough initial data $z_0+\langle\nabla\rangle^{-1}z_1\in H^s$. 
Left: $s=1/6$; \  right: $s=2/3$.\label{fig3}}
\end{figure}

\begin{figure}
\begin{center}
\subfigure{\includegraphics[width=0.49\textwidth]{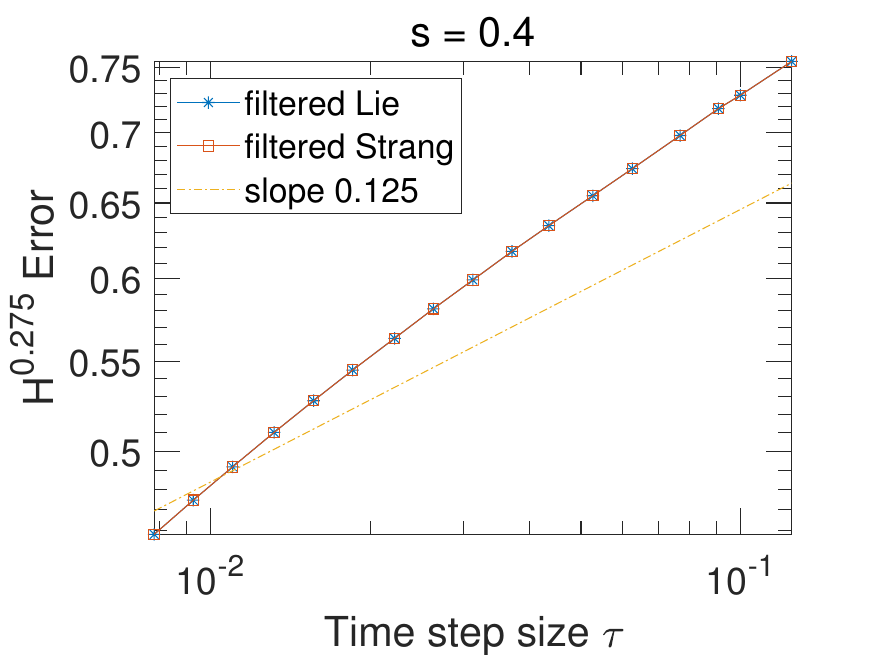}}
\subfigure{\includegraphics[width=0.49\textwidth]{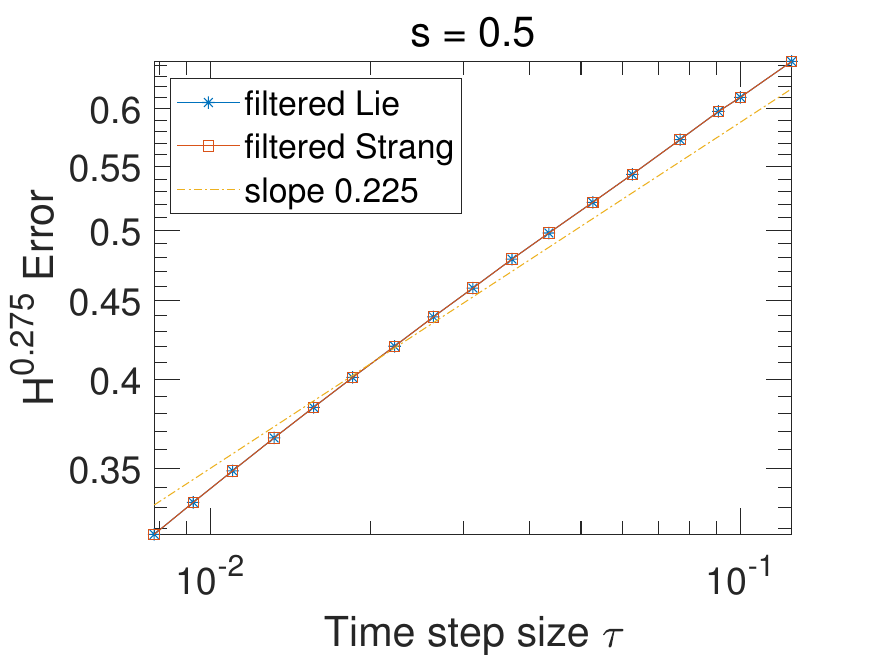}}
\end{center}
\caption{$H^\frac{11}{40}$ temporal error of the filtered Lie and Strang splitting schemes for the two-dimensional NKG with very rough initial data $z_0+\langle\nabla\rangle^{-1}z_1\in H^s$. 
Left: $s=0.4$; \ right: $s=0.5$.\label{fig4}}
\end{figure}

In Figure~\ref{fig3}, we performed numerical experiments for the Strang splitting scheme \eqref{kgstrang}. Here we normalized the initial data as $\|z_0+\langle\nabla\rangle^{-1}z_1\|_{L^2}=1$, and took Fourier modes $N=2^{14}$. For the reference solution, we used a small time step size $\tau=2^{-13}$. 
This is an experiment for the claim given in Remark~\ref{remsob}, where following the discussions of Section 3, the Strang splitting is supposed to have regular convergence behavior for $s>\frac12$ and irregular behavior for $s\in(\frac16,\frac12)$. From the figure on the left, we can see that the convergence behavior for the borderline case $s=\frac16$ is very irregular. However, with the gain of regularity of the initial data, we can see that the convergence behavior becomes regular again. For instance, when $s=\frac23>\frac12$, as we can see from the figure on the right, the convergence behavior is regular, and the convergence order is almost $2$.

In Figure~\ref{fig4}, we performed numerical experiments for the filtered Lie~\eqref{plie} and the filtered Strang splitting scheme \eqref{pstrang}. The initial data is normalized as $\|z_0+\langle\nabla\rangle^{-1}z_1\|_{H^\frac{11}{40}}=1$, and we took Fourier modes $N=(2^{13},2^{13})$. For the reference solution, we used a small time step size $\tau=2^{-12}$. 
From the figure, we can see that the convergence behaviors of the two methods are almost the same. In particular, the Strang splitting does not perform better than Lie. Moreover, the convergence order of the right figure (i.e., when $s=0.5$) is almost $0.225$, which fits Theorem~\ref{theobourg} very well. However, the figure on the left (i.e., $s=0.4$) shows that the convergence order is a little bit better than predicted by our theorem, the possible reason here is that the reference solution is not accurate enough since the method converges very slow~\cite{Jisiam,Jiima}.

{}
\end{document}